\documentclass[12pt]{amsart}
\usepackage{latexsym}
\usepackage{amsthm}
\usepackage{amssymb}
\usepackage{amsmath}
\usepackage{color}

\newtheorem{thm}{Theorem}[section]
\newtheorem{lem}[thm]{Lemma}
\newtheorem{prop}[thm]{Proposition}
\newtheorem{cor}[thm]{Corollary}

\theoremstyle{remark}
    \newtheorem*{rem}{Remark}
\theoremstyle{definition}
    \newtheorem{defn}[thm]{Definition}

\numberwithin{equation}{section}

\def\Z {{\mathbb Z}}
\def\N{{\mathcal N}}
\def\Q {{\mathbb Q}}
\def\T {{\tilde T}}
\def\D {{\mathcal D}}
\def\H {{\mathcal H}}
\def\h {{\mathfrak h}}
\def\l {\left\langle}
\def\r {\right\rangle}

\begin{document}

\title[centres of Hecke algebras]{On bases of centres of Iwahori--Hecke algebras of the symmetric group}
\author{Andrew Francis}
\address{School of Quantitative Methods and Mathematical Sciences, University of Western Sydney, NSW 1797, Australia}
\email[Andrew~Francis]{a.francis@uws.edu.au}
\author{Lenny Jones}
\address{Department of Mathematics, Shippensburg University, Pennsylvania, USA}
\email[Lenny~Jones]{lkjone@ship.edu}
\date{\today}
\begin{abstract}
In 1990, using norms, the second author constructed a basis for the centre of the Hecke algebra of the symmetric group $S_{n}$ over $\Q[\xi]$ \cite{Jon90}.
An integral ``minimal" basis was later given by the first author in 1999 \cite{Fmb}, following \cite{GR97}.
In principle one can then write elements of the norm basis as integral linear combinations of minimal basis elements.

In this paper we find an explicit non-recursive expression for the coefficients appearing in these linear combinations. These coefficients are expressed in terms of certain permutation characters of $S_n$.

In the process of establishing this main {theorem}, we prove the following {items} of independent interest: a result on the projection of the norms onto parabolic subalgebras, {the existence of} an inner product on the Hecke algebra with some interesting properties, and the existence of a partial ordering on the norms.
\end{abstract}

\maketitle

\setcounter{section}{-1}

\section{Introduction}

There are now three distinct descriptions of the centre of the Iwahori--Hecke algebra $\H$ of the symmetric group $S_n$.  It has two nice bases, one consisting of norms over $\Q[\xi]$ \cite{Jon90}, and {one} a ``minimal basis" of class elements over $\Z[\xi]$ \cite{GR97,Fmb}.  Thirdly, it is now known that the symmetric functions in Murphy elements are precisely the centre of $\H$ over {$\Z[\xi]$} \cite{Grah02}, and it follows that the elementary symmetric functions in Murphy elements generate the centre over $\Z[\xi]$.  A natural question is then to ask ``How are these {descriptions} related?".  The relationship between the elementary symmetric functions of Murphy elements and the minimal basis is now known precisely at least in one direction \cite{Fbaubles}, but relationships with the norm basis have been opaque.
Furthermore, the elucidation of the connections between the norm basis and the other bases is of interest since the norms of \cite{Jon90} are natural central structures which have been used to define Brauer-type homomorphisms for Hecke algebras \cite{Jon87,Du.GreenCor92,Fbrauer} and $q$--Schur algebras \cite{DD93}.

The goal of this paper is to describe an explicit relationship between the norm basis and the minimal basis for the centre of the Hecke algebra of the symmetric group $S_n$.  This relationship is given by an expression for the coefficients of class elements (the minimal basis) as they appear in the norms.  These coefficients are described in terms of the values of certain permutation characters of $S_n$.

Let $\alpha$ and $\lambda$ be partitions of $n$, with $w_\alpha$ an element of the conjugacy class $C_\alpha$ of $S_n$.  Let $l_\lambda$ and $l_\alpha$ be the lengths of the minimal elements in the corresponding conjugacy classes of $S_n$, and let $\xi$ be the defining indeterminate of the Hecke algebra.  Let $(1_{S_\lambda})^{S_n}$ be the permutation character of $S_n$ which arises from the induction to $S_n$ of the trivial character on the parabolic subgroup $S_\lambda$.  The main result is as follows.
\vspace{1mm}

\begin{minipage}{.9\textwidth}
\textbf{Theorem \ref{thm:coeff.Glambda.in.balpha}:}
Let $b_\alpha$ be an element of the norm basis and let $\Gamma_\lambda$ be an element of the minimal basis.  Then
\[
b_\alpha=\sum_{\lambda\vdash n}(1_{S_\lambda})^{S_n}(w_\alpha)\xi^{l_\lambda-l_\alpha}\Gamma_\lambda.\]
\end{minipage}
\vspace{2mm}

A considerable amount of machinery, which involves several results of independent interest, is developed in the course of obtaining Theorem \ref{thm:coeff.Glambda.in.balpha}.

The preliminary Section \ref{sec:defs} introduces most of the basic definitions and notation used throughout the paper.  The reader may wish to skim this section and return for reference as required later in the paper.
Section \ref{sec.double.cosets} contains results about double coset representatives of parabolic subgroups in the symmetric group $S_n$ which are required for Section \ref{sec:proj.norms.onto.maximal}.
A formula for the square of the Hecke algebra element corresponding to a distinguished double coset representative is given in Section \ref{sec:square.dcoset.Hecke}.
In {Section \ref{sec:bases.for.centre}, the main properties of the bases for the centre are briefly reprised}.
Section~\ref{sec:inner.prod} introduces an inner product on the Hecke algebra and gives some elementary properties.
In Section \ref{sec:norm.of.one}, we find the coefficient of $\Gamma_\lambda$ in $b_\alpha$ when $\alpha$ is trivial (Theorem \ref{thm:norm.of.1}), while in Section \ref{sec:partial.order.on.norms} we determine the coefficient of the Coxeter class element $\Gamma_{(n)}$ in $b_\alpha$ for all $\alpha\vdash n$ (Theorem \ref{thm:coeff.cox.elt.in.balpha}).
To establish Theorem \ref{thm:coeff.cox.elt.in.balpha}, we show that the basis of norms satisfies a partial order consistent with the refinement order on partitions (Theorem \ref{thm:balpha.partial.order}).
The descriptions of coefficients in Theorems \ref{thm:norm.of.1} and \ref{thm:coeff.cox.elt.in.balpha} are later made redundant by Theorem \ref{thm:coeff.Glambda.in.balpha}, but are necessary for {its proof}. 
Section \ref{sec:proj.norms.onto.maximal} gives the main projection theorem {(Theorem \ref{thm:proj.balpha.to.max.sgp})}, which uses a Mackey-type decomposition to give a rule for projecting norms onto a maximal parabolic subalgebra. This result has been used to study the Brauer homomorphism in \cite{Du.GreenCor92}.
In Section~\ref{sec:coeffs.of.class.elts}, Theorem \ref{thm:proj.balpha.to.max.sgp} is generalized to a rule for projecting onto arbitrary parabolic subalgebras {(Theorem \ref{thm:proj.balpha.trickle})}, and the main theorem quoted above is deduced.
Finally, the main result is demonstrated in Section \ref{sec:example} with some examples.

The authors thank the referee for many valuable comments and suggestions.

\section{Definitions and notation}\label{sec:defs}

Throughout we take $\mathbb N$ to mean the set of non-negative integers.

\subsection{Compositions, partitions and multipartitions}

A \emph{composition} $\lambda$ is {a finite ordered set} of positive integers.  If $\lambda=(\lambda_1,\dots,\lambda_r)$, the $\lambda_i$ are called the \emph{components} of $\lambda$.  If $\lambda$ is a composition we write $|\lambda|=\sum_{i=1}^r\lambda_i$.  If $|\lambda|=n$ we say $\lambda$ is a \emph{composition of $n$}, and we write $\lambda\vDash n$.
Two compositions are said to be \emph{conjugate} if they have the same components.

If $\lambda=(\lambda_1,\dots,\lambda_r)\vDash n$ then we define $\lambda-1$ to be the composition of $n$ obtained from $\lambda$ by replacing each $\lambda_i>1$ by the juxtaposed ordered pair of positive integers $\lambda_i-1$ and 1.  For example, if $\lambda=(3,4,1,7)$ then $\lambda-1=(2,1,3,1,1,6,1)$.

If $\lambda$ and $\mu$ are compositions of $n$ and either $\lambda=\mu$ or $\lambda$ can be obtained from $\mu$ by adding together adjacent components of $\mu$, we say $\mu$ is a \emph{refinement} of $\lambda$ and write $\mu\le\lambda$.

A \emph{partition} of $n$ is a composition whose components are weakly decreasing from left to right.  If $\lambda$ is a partition of $n$ we write $\lambda\vdash n$.

A \emph{multipartition} is a finite ordered set of partitions.
A \emph{$\lambda$--multipartition of $n$} for $\lambda=(\lambda_1,\dots,\lambda_r)\vDash n$ is an ordered set of partitions $\theta=(\theta_1,\dots,\theta_r)$ with $\theta_i\vdash\lambda_i$ for each $i=1,\dots, r$.
{Note that from any multipartition $\theta$ of $n$ we can derive a unique composition $\lambda$ of $n$ by removing the internal parentheses.
We call this unique composition $\lambda$, the \emph{derived composition} of the multipartition $\theta$.
For example, $\theta =((4,1),(3,2,1),(2,1))$ has derived composition $\lambda=(4,1,3,2,1,2,1)$ of $n=14$.}

By the \emph{components} of a $\lambda$--multipartition of $n$, we mean the components of its constituent partitions.
If $\alpha\vdash n$, then a \emph{$\lambda$--multipartition of $\alpha$} is a $\lambda$--multipartition of $n$ whose components are the components of $\alpha$.
Let $\Lambda_\lambda$ be the set of $\lambda$--multipartitions of $n$, and let $\Lambda_\lambda(\alpha)$ be the set of $\lambda$--multipartitions of $\alpha$.
For example, a $(3,5,2)$--multipartition of $(3,2,2,1,1,1)$ is $((2,1),(3,1,1),(2))$.

Note that for many choices of $\lambda$ and $\alpha$ there are no $\lambda$--multipartitions of $\alpha$; for instance there are no $(3,2)$--multipartitions of $(4,1)$.

\subsection{The symmetric group}\label{subsec:symm.group}

Let $S_n$ be the symmetric group on $n$ letters with generating set of simple reflections
\[S:=\{s_i=(i\ \ i+1)\mid 1\le i\le n-1\}.\]
{We use} both the $s_i$ notation and the cycle notation as expedient.  We adopt the convention that $S_{0}=S_{1}=\left\{ (1)\right\}$.

We say an expression for $w\in S_n$ is \emph{reduced} if there is no way to write $w$ as a word in fewer generators.  In this case we say the length $\ell(w)$ of $w$ is this minimal number of generators.
Symmetric groups act on sets of vectors in Euclidean space known as \emph{root systems}.  One can define the concepts of \emph{positive} and \emph{negative} roots which in turn can be used to describe the length of an element of $S_n$.  In particular if $\Phi^+$ and $\Phi^-:=\{-v\mid v\in\Phi^+\}$ are the sets of positive and negative roots respectively, then $\ell(w)=|w(\Phi^+)\cap\Phi^-|$.  A set of positive roots for the root system of $S_n$ is the set $\Phi^+=\{e_i-e_j\mid 1\le i<j\le n\}$ where $\{e_i\mid 1\le i\le n\}$ is the set of standard basis vectors for $\mathbb R^n$.

For $\lambda=(\lambda_1,\dots,\lambda_k)\vDash n$ define
\[S_\lambda:=S_{\lambda_1}\times S_{\lambda_2}\times\dots\times S_{\lambda_k},\]
where for $\lambda_i>1$, $S_{\lambda_i}$ is the subgroup of $S_n$ generated by the set $\{s_{\lambda_1+\dots+\lambda_{i-1}+1},\dots,s_{\lambda_1+\dots+\lambda_i-1}\}$, and for $\lambda_i=1$, $S_{\lambda_i}$ is the trivial subgroup.  Such a subgroup $S_\lambda$ is called a \emph{parabolic} subgroup of $S_n$.  Note that $S_\mu\le S_\lambda$ if and only if $\mu\le\lambda$ (that is, $\mu$ is a refinement of $\lambda$).

If $\theta=(\theta_1,\dots,\theta_t)$ is a multipartition (the $\theta_i$ are partitions), then set $S_\theta:=S_{\theta_1}\times\dots\times S_{\theta_t}$.

The unique element of a $S_\lambda$--$S_\mu$ double coset of $S_n$ of minimal length is called a \emph{distinguished} double coset representative (such elements are well-known to be unique --- see \cite{Carter85}). Let $\D_{\lambda\mu}$ denote the set of distinguished $S_\lambda$--$S_\mu$ double coset representatives in $S_n$.

The conjugacy classes of $S_n$ are indexed by partitions $\lambda$ of $n$.  Write $C_\lambda$ for the conjugacy class consisting of elements of $S_n$ of cycle type $\lambda$.  Write $\underline{C_\lambda}$ for the sum of elements in the conjugacy class $C_\lambda$.
If $\lambda=(\lambda_1,\dots,\lambda_r)$ then set
\[w_\lambda=(s_1\dots s_{\lambda_1-1})(s_{\lambda_1+1}\dots s_{\lambda_1+\lambda_2-1})\dots(s_{\lambda_1+\dots+\lambda_{r-1}+1}\dots s_{\lambda_1+\dots+\lambda_{r}-1})
\]
where we take each empty sequence of $s_i$'s (when $\lambda_j=1$) to be the identity.
Then $w_\lambda$ is a \emph{Coxeter element} of the subgroup $S_\lambda$, and also a minimal length element of the conjugacy class $C_\lambda$ in $S_n$.

In $S_\lambda$ for $\lambda\vDash n$, the conjugacy classes are indexed by the set $\Lambda_\lambda$ of $\lambda$--multipartitions of $n$.
In particular, in $S_{(k,n-k)}$, the classes are indexed by $(k,n-k)$--multipartitions of $n$.
If $\lambda\vDash n$ and $\theta\in\Lambda_\lambda$, then $C_\theta$ denotes the conjugacy class in $S_\lambda$ corresponding to the composition of $n$ {derived from} $\theta$.

As usual, $|C_\theta|_{S_\lambda}$ denotes the size of the conjugacy class $C_\theta$ in $S_\lambda$.

Let $l_\lambda$ be the length of a shortest element of the conjugacy class $C_\lambda$, that is, $l_\lambda=\ell(w_\lambda)$.

For $w\in S_{n}$ and {fixed} $k\in \{1,\dots,n-1\}$, define {$\#(w)$} to be the minimal number of times the generator $s_k=(k\ \ k+1)$ must appear in any reduced expression for $w$. {Unless otherwise noted, $k$ is assumed to be fixed throughout this paper.} For $w\in S_n$ and $I\subseteq\{1,\dots,n\}$, we write $w.I$ for the image of the action of $w$ on the set $I$.

The \emph{Bruhat order} on $S_n$ is defined as follows.  For $v,w\in S_n$, we say  $v\le w$ if there exists a reduced expression of $v$ which is a subword of a reduced expression for $w$.

For any groups $H\le G$, we use the standard notation $C_{G}(H)$
and $N_{G}(H)$ to indicate the \emph{centralizer} and
\emph{normalizer} respectively, of $H$ in $G$.

\subsection{The Hecke algebra}\label{subsec:Hecke.algebras}

In this paper {we use} the normalized version
of the generators for the Hecke algebra, giving us an algebra over the ring $\Z[\xi]$, where $\xi$ is an indeterminate.  The exact connection between this definition and the standard definition over $\Z[q^{1/2},q^{-1/2}]$ is given {in the remark} below.

The Iwahori--Hecke algebra $\H:=\H_n$ of $S_n$ is the associative $\Z[\xi]$-algebra generated by the set $\{\T_s\mid s\in S\}$ with identity $\T_1$ and subject to the relations
\begin{align*}
 \T_{s_i}\T_{s_j}&=\T_{s_j}\T_{s_i} &\text{if $|i-j|\ge 2$}\ \\
 \T_{s_i}\T_{s_{i+1}}\T_{s_i}&=\T_{s_{i+1}}\T_{s_i}\T_{s_{i+1}} &\text{for $1\le i\le n-2$}\ \\
 \T_{s_i}^2&=\T_1+\xi\T_{s_i}   &\text{for $1\le i\le n-1$.}
\end{align*}
If $w=s_{i_1}\dots s_{i_r}$ is a reduced expression for $w$, then we write $\T_w:=\T_{s_{i_1}}\dots \T_{s_{i_r}}$.  Then $\H$ is a free $\Z[\xi]$--module with basis $\{\T_w\mid w\in S_n\}$.

When specialization of $\H$ to $\xi=0$ is used in this paper it will be assumed that the specialization is to the  group algebra $\Z S_n$.  If $h\in\H$ we write $h|_{\xi=0}$ for the specialization of $h$ at $\xi=0$.

If $\lambda \vDash n$, we let $\H_{\lambda}$ denote the parabolic subalgebra of $\H_{n}$ generated by $\{ \T_{s} \mid s\in S^{\prime}\}$, where $S^{\prime}$ is the subset of $S$ consisting of the simple reflections which generate the parabolic subgroup $S_{\lambda}$. {For any multipartition $\theta$ of $n$ with derived composition $\lambda$ of $n$, we define $\H_{\theta}:=\H_{\lambda}$.}
\begin{rem} Set $\T_s:=q^{-1/2}T_s$ for $s\in S$ and $\xi=q^{1/2}-q^{-1/2}$.  Then $\H$ is a subalgebra of $\H_q$, the more standard Hecke algebra generated by $\{T_s\mid s\in S\}$ over $\Z[q^{1/2},q^{-1/2}]$.  A principal reason for defining the algebra with normalized generators is that doing so gives $\H$ a natural positivity and an associated partial order on the positive cone.  Many results on the centre of $\H$ have more natural statements and proofs when the algebra is defined in this way.  The main results of this paper are all readily translated back to statements over $\Z[q^{1/2},q^{-1/2}]$.
\end{rem}

Let $\H^+=\sum_{w\in S_n}{\mathbb N}[\xi]\T_w$.  For $a,b\in\H^+$ we say $a\le b$ when $b-a\in\H^+$.  {If in addition $b-a\neq 0$ then we write $a<b$.}
The partial order restricts to the positive cone of the base ring, $\mathbb N[\xi]=\Z[\xi]^+$.  For $a,b\in\mathbb N[\xi]$, we say $a\le b$ when $b-a\in\mathbb N[\xi]$, and $a<b$ if in addition $b-a\neq 0$.

If $A$ is a subalgebra of $\H$, then the \emph{centralizer of $A$ in $\H$}, denoted $Z_\H(A)$, is the set of elements in $\H$ which commute with every element of $A$.  The \emph{centre} of $\H$ is $Z(\H):=Z_\H(\H)$.  Set  $Z(\H)^+:=Z(\H)\cap\H^+$.

We say that an element $\T_{w}$ in $\H_{n}$ {\emph{contains}}
a particular generator $s\in S$ when $s\le w$ in the Bruhat order.
We also say (with some abuse of language) that $h=\sum_{w\in S_n}r_w\T_w\in\H_n$ \emph{contains} $\T_w$, or that $\T_w$ \emph{occurs} in $h$, if $r_w\neq 0$.

\section{Double cosets of maximal parabolic subgroups}\label{sec.double.cosets}
{This section and the next contain a number of results needed for later sections, most of which appear in \cite{Jon87} but not, as far as we can see, in the available literature.}

{Recall that $k\in \{1,\dots,n-1\}$ is fixed. Also, throughout this section we let $M=$ min$\left\{ k,n-k \right\}$ and $\lambda=(k,n-k)\vDash n$.
For $0\le m\le M$ and $w\in S_n$, set
 \begin{align*}
 d_m&:=(k\ \ k+m)(k-1\ \ k+m-1)\cdots (k+1-m\ \ k+1),\\
 D&:=\{d_m\mid 0\le m\le M\},\\
||w||&:=|w.\{1,\dots,k\}\cap \{k+1,\dots,n\}|.
\end{align*}}

\begin{prop}
\label{prop:ddcosets.reps.3.2}
{\text{}
\begin{enumerate}
  \item The elements $d_m$ satisfy the following:
\begin{enumerate}
\item\label{item:ddcoset.reps.are.invols} $d_{m}^{-1}=d_{m}$
\item\label{item:ddcoset.reps.length} $\ell (d_{m})=m^{2}$
\item\label{item:ddcoset.reps.hash} $\# (d_{m})=m=||d_{m}||$.
\end{enumerate}
  \item\label{item:ddcoset.reps} $D=\D_{\lambda\lambda}$.
  \end{enumerate}}
\end{prop}

\begin{proof}
Part \eqref{item:ddcoset.reps.are.invols} is obvious since $d_{m}$ is a product of disjoint transpositions.

The results of \eqref{item:ddcoset.reps.length} and \eqref{item:ddcoset.reps.hash} are trivial for $m=0$ and $m=1$, so we assume $m\ge 2.$
Using a standard root system argument (see for example \cite{Carter85}), observe that $d_{m}$ takes the
set $\{e_{k-i}-e_{k+1+j}\mid 0\le i,j\le m-1\}$ of $m^{2}$ positive roots to negative roots, which implies that $\ell(d_{m})\ge m^{2}$.
Since $d_m=d_{m-1}\left(s_{k-m+1}\dots s_{k-1}s_{k+m}\dots s_k\right)$, we have that
$\ell (d_{m})\le \ell (d_{m-1})+2m-1$. By induction, $\ell (d_{m-1})=(m-1)^{2}.$ Hence $\ell (d_{m})\le (m-1)^{2}+2m-1=m^{2}$, giving $\ell(d_{m})=m^{2}$, which proves part \eqref{item:ddcoset.reps.length}.

As above, write $d_m=d_{m-1}w s_k$ where $\# (w)=0.$ By induction assume that $\# (d_{m-1})=m-1.$ Therefore, $\# (d_m)\le m$. {From the definition of $d_m$ it is immediate that $||d_m||=m.$ Then, since $||d_m||\le \# (d_m)$, we conclude that $\# (d_m)=m=||d_m||$ as required for \eqref{item:ddcoset.reps.hash}.}

Note that if $w$ and $u$ are in the same $S_\lambda$-$S_\lambda$ double coset of $S_n$, then $||w||=||u||$.  As $||d_m||=m$ for each $m$, we have that each $d_m$ lies in a distinct $S_\lambda$-$S_\lambda$ double coset.  Observe that {$|D|=M+1$}, and since there are {$M+1$} distinct $S_\lambda$-$S_\lambda$ double cosets (\cite[Theorem 1.3.10]{JK81}), $D$ is a set of double coset representatives.
Since any element $w$ of $S_\lambda d_mS_\lambda$ takes the $m^2$ above--mentioned positive roots to negative roots, {it follows that} $\ell(w)\ge m^2$. Thus $d_m$ is distinguished, proving \eqref{item:ddcoset.reps}.
\end{proof}

\begin{cor}
\label{cor:x.in.Xm.iff.hash.3.4}
{Let $x\in S_{n}$. Then
\begin{enumerate}
\item $\#(x)=||x||$
\item $x\in S_{\lambda}d_mS_{\lambda}$ if and only if $\# (x)=m.$
\end{enumerate}}
\end{cor}
\begin{proof}
{Suppose that $x\in S_{\lambda}d_mS_{\lambda}$. Then we can write $x=yd_{m}z$ with $\# (y)=0=\# (z)$, and clearly, $\# (yd_{m}z)\le \#(d_m)$. Recall that $\#(d_m)=||d_m||$ from Proposition \ref{prop:ddcosets.reps.3.2}(1c). Since $x$ and $d_m$ are both elements of $S_{\lambda}d_mS_{\lambda}$, we have $||x||=||d_m||$. Hence, $\# (x)=\# (yd_{m}z)\le \#(d_m)=||d_m||=||x||$. Note that each occurrence of $s_k$ in any expression for $x$ produces at most one element of $x.\{1,\dots,k\}\cap \{k+1,\dots,n\}$. Thus, $\#(x)\ge ||x||$, which proves (1).}

{We use part (1) and arguments in its proof to establish (2). For $x\in S_{\lambda}d_mS_{\lambda}$, we have $\#(x)=||x||=||d_m||=m$.
Conversely, if $\# (x)=m,$ then $||x||=m=||d_m||$, which implies that $x\in S_{\lambda}d_mS_{\lambda}$, and the proof is complete.}
\end{proof}

For {$0\le m\le M$}, let
\[P_m:=S_{\lambda}\cap S_{\lambda}^{d_m}.\]

\begin{prop}
    \label{prop:Pm=etc.3.7}\text{}
\begin{enumerate}
  \item $P_m=P_m^{d_m}$.\label{cor:Pm.to.dm=Pm}
  \item $P_m=S_{(k-m,m,m,n-k-m)}$.\label{prop:part:Pm=S..etc}
\end{enumerate}
\end{prop}

\begin{proof}
Part \eqref{cor:Pm.to.dm=Pm} is immediate from Proposition \ref{prop:ddcosets.reps.3.2}{\eqref{item:ddcoset.reps.are.invols}}.

Clearly $S_{(k-m,m,m,n-k-m)}\le S_{\lambda}$.  Observe that $d_{m}$ normalizes $S_{(k-m,m,m,n-k-m)}$, and so
\[
S_{(k-m,m,m,n-k-m)}=S_{(k-m,m,m,n-k-m)}^{d_m}\le S_{\lambda}^{d_m},
\]
giving $S_{(k-m,m,m,n-k-m)}\le P_m$.

If $s_{k-m}\in S_{\lambda}^{d_m}$, then $d_ms_{k-m}d_m\in S_{\lambda}$.
But
\[
d_ms_{k-m}d_{m}=s_{k-m}s_{k-m+1}\cdots s_{k-1}s_ks_{k-1}\cdots s_{k-m+1}s_{k-m}
\]
 which implies that $s_k\in S_{\lambda}$, a contradiction. Hence, $s_{k-m}\not \in S_{\lambda}^{d_m}$.  Similarly, $s_{k+m}\not \in S_{\lambda}^{d_m}$.  Since $P_m$ is parabolic
it follows that $P_m\le S_{(k-m,m,m,n-k-m)}$, and \eqref{prop:part:Pm=S..etc} is proved.
\end{proof}

\begin{cor}
\label{cor:dwd.is.simple.refl.3.8}\label{cor:length.dxd=length.x.3.9}
If $w\in P_{m}$ then $\ell(d_{m}wd_{m})=l(w)$.
\end{cor}
\begin{proof} Immediate from Proposition \ref{prop:Pm=etc.3.7}\eqref{prop:part:Pm=S..etc} and the definition of $d_m$.
\end{proof}

\begin{cor}
\label{cor:drd.is.coset.rep}
Let $U$ be a parabolic subgroup of $P_{m}$.
Let $r$ be a distinguished right coset representative for $U$ in $P_{m}$.
Then $d_{m}rd_{m}$ is a distinguished right coset representative for $U^{d_{m}}$ in $P_{m}.$
\end{cor}
\begin{proof} This is immediate from Corollary \ref{cor:length.dxd=length.x.3.9}.
\end{proof}

\section{The square of the Hecke algebra element corresponding to a distinguished double coset representative}\label{sec:square.dcoset.Hecke}\text{}

The goal of this section is to prove Proposition \ref{prop:Tdm.squared}, which gives an expansion for the square of the Hecke algebra element corresponding to a distinguished double coset representative of a maximal parabolic subgroup.  Proposition \ref{prop:Tdm.squared} forms part of the machinery needed for our analysis in Section \ref{sec:proj.norms.onto.maximal} of the projection of the norm basis onto a maximal parabolic subalgebra. {Throughout this section we let $\lambda=(k,n-k)\vDash n$.}

For $x,y,z\in S_n$, define the polynomial $f_{xyz}\in {\mathbb N}[\xi]$ to be the coefficient of $\T_z$ occurring in the expansion of $\T_x\T_y$.  That is, write
\begin{equation}\T_x\T_y=\sum_{z\in S_n}f_{xyz}\T_z.\label{eq:TxTy=SumfxyzTz}\end{equation}

\begin{lem}
  Let $x,y,z\in S_n$, with $f_{xyz}\neq 0$.
  \begin{enumerate}
    \item $\#(x)-1\le\#(s_k x)\le\#(x)+1$.
    \item If $\#(x)\ne\#(y)$ then $\#(z)\ge 1$.
  \end{enumerate}
\end{lem}

\begin{proof}
Firstly note that (2) follows immediately from (1) by induction on the length of $x$ or $y$.  For (1), the upper bound is clear.

We will use induction on $\#(x)=t$ to establish the lower bound in (1).  Write $x=w_1d_tw_2$, with $w_1,w_2\in S_{\lambda}$ (so that $\#(w_1)=\#(w_2)=0$), and assume that $\# (s_kw_1d_tw_2)<t-1$.  Then $s_kw_1d_tw_2\in S_{\lambda}d_vS_{\lambda}$ for some  $v<t-1.$  That is, $s_kw_1d_tw_2=\hat{x}$ for some $\hat{x}\in S_{\lambda}d_vS_{\lambda}.$ Then $d_{t}=w_1^{-1}s_k\hat{x}w_2^{-1}.$
By induction,
$\# (\hat{x})-1 \le \# (s_k\hat{x})\le \# (\hat{x}) +1.$ In other words, $v-1\le \# (s_k\hat{x})\le v+1.$ Since $w_1^{-1},w_2^{-1}\in S_{\lambda}$ we have that $d_{t}$ is an element of either $S_{\lambda}d_{v-1}S_{\lambda}$, $S_{\lambda}d_{v}S_{\lambda}$ or $S_{\lambda}d_{v+1}S_{\lambda}$, each of which is impossible since $v<t-1.$
\end{proof}

We recall a result of Shi \cite{Shi.Bruhat90}.

\begin{thm}[\protect{\cite[Theorem 8]{Shi.Bruhat90}}]\label{thm.Shi.Bruhat}
  Let $x,y,z\in S_n$. If $f_{xyz}\neq 0$ then $xy\le z$ in the Bruhat order.
\end{thm}

An elementary consequence is the following:
\begin{cor}\label{lem:sharp.in.TxTy}\label{cor:coeff.tau.3.13}\label{lem:tau.in.prod.with.Bruhat}\label{lem:coeff.tau.in.H.with.length.3.15}
Let $x,y,z \in S_{n}$.  If $f_{xyz}\neq 0$, then
\begin{enumerate}
  \item $\# (xy)\le \# (z)$.
  \item If $\# (x)=t$ and $\# (y)=0$ then $\# (z)=t.$
\end{enumerate}
\end{cor}

\begin{lem}\label{lem:slant.product}
For $j_1\le j_2$,
\[\T_{s_{j_1}\dots s_{j_2}}\T_{s_{j_2}\dots s_{j_1}}=\T_1+\xi\sum_{i=j_1}^{j_2}\T_{s_{j_1}\dots s_i\dots s_{j_1}}.\]
\end{lem}
\begin{proof}
  Elementary (by induction on $j_2$).
\end{proof}

\begin{prop}\label{prop:Tdm.squared}
\[
(\T_{d_{m}})^{2}=\T_{1}+ \sum_{w \in S_{n}}f_{w}\T_{w}
\]
where $\# (w)\ge 1$ for $0\neq f_{w}\in{\mathbb N}[\xi]$.

(Note that we abbreviate $f_{d_md_mw}$ in the expansion of \eqref{eq:TxTy=SumfxyzTz} to $f_w$.)
\end{prop}
\begin{proof}
For $m=1$, the proposition is easily verified. Assume $m\ge 2$, and write $d_{m}=d_{m-1}uvs_k$ where $u=s_{k-m+1}s_{k-m+2}\cdots s_{k-1}$
and $v=s_{k+m-1}s_{k+m-2}\cdots s_{k+1}$.  Then, since $d_{m}^{-1}=d_{m}$ and $uv=vu$, we have
\[
\begin{array}{rl}
(\T_{d_m})^2&=\T_{d_{m-1}}\T_{u}\T_{v}\T_{s_{k}}
\T_{s_k}\T_{v^{-1}}\T_{u^{-1}}\T_{d_{m-1}}\\
&=\T_{d_{m-1}}(\T_{u}\T_{u^{-1}}\T_{v}\T_{v^{-1}})\T_{d_{m-1}}+\xi \T_{d_{m}}\T_{v^{-1}u^{-1}d_{m-1}}.
\end{array}
\]
Now, $v^{-1}u^{-1}d_{m-1}\in S_{\lambda}d_{m-1}S_{\lambda}$ and so by Corollary \ref{cor:x.in.Xm.iff.hash.3.4}{(2)}, $\# (v^{-1}u^{-1}d_{m-1})=m-1.$ Since $\# (d_{m})=\# (d_{m-1})+1$, Corollary \ref{cor:coeff.tau.3.13} implies that every term of the product $\T_{d_{m}}\T_{v^{-1}u^{-1}d_{m-1}}$ {contains ${s_k}$}.  Now $\T_{u}\T_{u^{-1}}=\T_1+\xi\sum_{i=k-m+1}^{k-1}\T_{s_{k-m+1}\dots s_i\dots s_{k-m+1}}$ by Lemma \ref{lem:slant.product}.
Similarly, $\T_{v}\T_{v^{-1}}=\T_{1}+\xi\sum_{i=2}^{m}\T_{s_{k+i-1}\dots s_{k+m-1}\dots s_{k+i-1}}$.
Hence
\begin{multline*}
\T_{d_{m-1}}(\T_{u}\T_{u^{-1}}\T_{v}\T_{v^{-1}})\T_{d_{m-1}}\\
=\T_{d_{m-1}}\left(\T_1+\xi\sum_{i=k-m+1}^{k-1}\T_{s_{k-m+1}\dots s_i\dots s_{k-m+1}}\right)\\
\times\left(\T_{1}+\xi\sum_{i=2}^{m}\T_{s_{k+i-1}\dots s_{k+m-1}\dots s_{k+i-1}}\right)\T_{d_{m-1}}.
\end{multline*}

By induction, assume
\[
    (\T_{d_{m-1}})^{2}=\T_{1}+ \sum_{w \in S_{n}}f'_{w}\T_{w}
\]
where $\# (w)\ge 1$ for $f'_{w}\ne 0$, and again we have abbreviated the coefficient $f_{d_{m-1}d_{m-1}w}$ from \eqref{eq:TxTy=SumfxyzTz} to $f'_w$.

The proposition will be proved if we can show that all terms in each of the products
\[
    \T_{d_{m-1}}\T_{s_{k-m+1}\dots s_{k-m+a-1}\dots s_{k-m+1}}\T_{d_{m-1}},\quad
    \T_{d_{m-1}}\T_{s_{k+b-1}\dots s_{k+m-1}\dots s_{k+b-1}}\T_{d_{m-1}}
\]
and
\[
    \T_{d_{m-1}}\T_{s_{k-m+1}\dots s_{k-m+a-1}\dots s_{k-m+1}}\T_{s_{k+b-1}\dots s_{k+m-1}\dots s_{k+b-1}}\T_{d_{m-1}},
\]
contain $s_k$, for $2\le a,b\le m$.

The product $\left(s_{k-m+1}\dots s_{k-m+a-1}\dots s_{k-m+1}\right)\left(s_{k+b-1}\dots s_{k+m-1}\dots s_{k+b-1}\right)$ adds in length, and since it is an element of $S_{\lambda}$ for all $a$ and $b$, its product with $d_{m-1}$ is length-additive.  Hence,
\begin{multline*}
 \T_{s_{k-m+1}\dots s_{k-m+a-1}\dots s_{k-m+1}}\T_{s_{k+b-1}\dots s_{k+m-1}\dots s_{k+b-1}}\T_{d_{m-1}}\\
=\T_{s_{k-m+1}\dots s_{k-m+a-1}\dots s_{k-m+1}s_{k+b-1}\dots s_{k+m-1}\dots s_{k+b-1}d_{m-1}}.
\end{multline*}
It is easily determined that
\begin{align*}
d_{m-1}s_{k-m+1}\dots s_{k-m+a-1}\dots s_{k-m+1}d_{m-1}&\in S_{\lambda}d_{1}S_{\lambda},\\
d_{m-1}s_{k+b-1}\dots s_{k+m-1}\dots s_{k+b-1}d_{m-1}&\in S_{\lambda}d_{1}S_{\lambda}
\end{align*}
and
\begin{multline*}
   d_{m-1}(s_{k-m+1}\dots s_{k-m+a-1}\dots s_{k-m+1})\cdot \\
   (s_{k+b-1}\dots s_{k+m-1}\dots s_{k+b-1})d_{m-1}\in S_{\lambda}d_{2}S_{\lambda}
\end{multline*}

for any $a$ and $b$ with $m\ge 2.$  That is, by Corollary \ref{cor:x.in.Xm.iff.hash.3.4}{(2)}
\begin{multline*}
    \#(d_{m-1}s_{k-m+1}\dots s_{k-m+a-1}\dots s_{k-m+1}d_{m-1})=1\\
    =\# (d_{m-1}s_{k+b-1}\dots s_{k+m-1}\dots s_{k+b-1}d_{m-1}
\end{multline*}
and
\begin{multline*}
    \#\big(d_{m-1}(s_{k-m+1}\dots s_{k-m+a-1}\dots s_{k-m+1})\\ \cdot(s_{k+b-1}\dots s_{k+m-1}\dots s_{k+b-1})d_{m-1}\big)=2.
\end{multline*}
Thus an application of Corollary \ref{lem:coeff.tau.in.H.with.length.3.15}{(1)} to each of the products
\[
\T_{d_{m-1}}\T_{s_{k-m+1}\dots s_{k-m+a-1}\dots s_{k-m+1}d_{m-1}},\quad
\T_{d_{m-1}}\T_{s_{k+b-1}\dots s_{k+m-1}\dots s_{k+b-1}d_{m-1}}
\]
{and}
\[
\T_{d_{m-1}}\T_{(s_{k-m+1}\dots s_{k-m+a-1}\dots s_{k-m+1})(s_{k+b-1}\dots s_{k+m-1}\dots s_{k+b-1})d_{m-1}}
\]
completes the proof.
\end{proof}

\begin{cor}\label{cor:tau.in.Tdm.Tgam.Tdm.3.17}
Let $v\in P_{m}.$ Then
\[
    \T_{d_m}\T_{v}\T_{d_m}=\T_{d_{m}v d_{m}}+ \sum_{w \in S_{n}}g_{w}\T_{w}
\]
where $\# (w)\ge 1$ for $0\neq g_{w}\in\mathbb N[\xi]$.
\end{cor}
\begin{proof} Since $v\in S_{\lambda}$,
we have $\ell (v d_{m})=\ell (v) +\ell(d_{m})$.  From Proposition \ref{prop:Pm=etc.3.7}\eqref{cor:Pm.to.dm=Pm}, we have that $\hat{v}=d_{m}v d_{m}\in P_{m}$, so also $\ell (d_{m}\hat{v})=\ell (d_{m})+\ell (\hat{v}).$
Thus
\begin{align*}
\T_{d_{m}}\T_{v}\T_{d_{m}}
&=\T_{d_{m}}\T_{d_{m}}\T_{\hat{v}}\\
&=(\T_{1}+ \sum_{u\in S_{n}}f_u\T_u)\T_{\hat{v}}\quad \text{(where $\# (u)\ge 1$ for $f_u\ne 0$, by}\\
&\hspace{3cm}\text{Proposition \ref{prop:Tdm.squared})}\\
&= \T_{\hat{v}}+ \sum_{u \in S_n}f_u\T_u\T_{\hat{v}}\\
&= \T_{d_{m}v d_{m}}+ \sum_{w \in S_{n}}g_{w}\T_{w}
\end{align*}
where $\# (w)\ge 1$ when $g_{w}\ne 0$, by Corollary \ref{lem:tau.in.prod.with.Bruhat}{(2)} since $\# (\hat{v})=0.$
\end{proof}

\section{Bases for the centre of the Hecke algebra}\label{sec:bases.for.centre}

In this section we introduce the two bases for the centre of the Hecke algebra whose relationship is the main topic of this paper.
With the exception of Proposition \ref{lem:norms.from.bips}, the results in this section appear in the existing literature exactly as stated here, or in slightly less generality.

Some results from \cite{Jon90} have been restated in the context of the Hecke algebra over $\Z[\xi]$, and in the generality of compositions rather than partitions where appropriate.

\subsection{The norm basis}\label{subsec:norm.basis}

\begin{defn}
For $h\in\H$, $\lambda,\mu,\alpha\vDash n$ and $\mu\le\lambda$, we define the \emph{relative norm} of $h$ from $S_\mu$ to $S_\lambda$ to be
\[\N_{S_\lambda,S_\mu}(h):=\sum_{d\in\D}\T_{d^{-1}}h\T_d,\]
where $\D$ is the set of distinguished right coset representatives of $S_\mu$ in $S_\lambda$.
In addition, define
\[\eta_\lambda:=\N_{S_{\lambda-1},1}(\T_{w_\lambda})\]
and
\[b_\alpha:=\N_{S_n,S_\alpha}(\eta_\alpha).\]
\end{defn}
{For any multipartition $\theta$ with derived composition $\lambda$, we define $\eta_\theta:=\eta_\lambda$.}

As $C_{(n)}$ is the Coxeter class of $S_n$, we call $\eta_{(n)}$ the Coxeter class element of $\H_n$.  Similarly, $\eta_\lambda$ is the Coxeter class element of $\H_\lambda$.

\begin{prop}\label{lem:norms.from.bips}
If $\alpha$ and $\beta$ are conjugate compositions, then $b_\alpha=b_\beta$.
\end{prop}
\begin{proof}
Since any pair of conjugate compositions can be obtained from one another via
a sequence of exchanges of adjacent components, it suffices to consider two conjugate compositions which differ by a single adjacent pair.  That is, we may assume $\alpha=(\lambda_1,\dots,\lambda_k,\lambda_{k+1},\dots,\lambda_r)$ and $\beta=(\lambda_1,\dots,\lambda_{k+1},\lambda_k,\dots,\lambda_r)$.  Then
\begin{equation*}
b_\alpha=\N_{S_n,S_\alpha}(\eta_\alpha)
=\N_{S_n,S_{(\lambda_1,\dots,\lambda_k+\lambda_{k+1},\dots,\lambda_r)}}\left(\N_{S_{(\lambda_1,\dots,\lambda_k+\lambda_{k+1},\dots,\lambda_r)},S_\alpha}(\eta_\alpha)\right),
\end{equation*}
with
{
\begin{multline*}
\N_{S_{(\lambda_1,\dots,\lambda_k+\lambda_{k+1},\dots,\lambda_r)},S_\alpha}(\eta_\alpha)\\ =\eta_{(\lambda_1,\dots,\lambda_{k-1})}\N_{S_{\lambda_k+\lambda_{k+1}},S_{(\lambda_{k},\lambda_{k+1})}}(\eta_{(\lambda_k,\lambda_{k+1})})\eta_{(\lambda_{k+2},\dots,\lambda_r)}. \end{multline*}    }
So to show $b_\alpha=b_\beta$ it suffices to show
\begin{equation}\label{eq:norms.from.bips}
\N_{S_{\lambda_k+\lambda_{k+1}},S_{(\lambda_{k},\lambda_{k+1})}}(\eta_{(\lambda_k,\lambda_{k+1})})=
\N_{S_{\lambda_k+\lambda_{k+1}},S_{(\lambda_{k+1},\lambda_{k})}}(\eta_{(\lambda_{k+1},\lambda_{k})}).
\end{equation}
These norms are both central in $\H_{\lambda_{k}+\lambda_{k+1}}$ by \cite[Proposition 2.13]{Jon90}, and they are images of each other under the algebra automorphism of $\H_{\lambda_{k}+\lambda_{k+1}}$ defined by reflecting the Dynkin diagram about its midpoint.  This automorphism fixes central elements and hence \eqref{eq:norms.from.bips} holds, proving the lemma.
\end{proof}

{Note that if $\alpha$ and $\beta$ are conjugate compositions of $n$, there exists $x\in S_{n}$ such that $S_\alpha^x=S_\beta.$ Therefore, we write $\alpha^x=\beta$ in this situation. As before, {for any multipartition $\theta$ of $n$ with derived} composition $\lambda$ of $n$, we define $\theta^x:=\lambda^x.$}
\begin{thm}[\cite{Jon90}]\label{lem:etacox}
Let {$\lambda\vDash n$}.
 \begin{enumerate}
  \item $\eta_\lambda\in Z(\H_\lambda)$.
  \item Every term with non-zero coefficient in the element $\eta_{(n)}$ contains every simple reflection in its reduced form.
  \item $\eta_{(n)}$ contains terms of length $n-1$ and greater, and specializes on $\xi=0$ to the Coxeter class sum in $\Z S_n$.
 \end{enumerate}
\end{thm}
\begin{proof}
Part (1) is \cite[Lemma (3.23)(3)]{Jon90}.
Part (2) is \cite[Lemma (3.14)]{Jon90}.
Part (3) is \cite[Corollary (3.15) and Lemma (3.23)(4)]{Jon90}.
\end{proof}

\begin{thm}[\cite{Jon90}]\label{thm:Jones.norm.properties}
Let {$\lambda\vDash n$}.
  \begin{enumerate}
  \item (Transitivity) If $S_\mu\le S_\lambda$ then $\N_{S_n,S_\mu}(h)=\N_{S_n,S_\lambda}\left(\N_{S_\lambda,S_\mu}(h)\right)$.\label{thm:part:norm.transitivity}
  \item If $h\in Z(\H_\lambda)$ then $\N_{S_n,S_\lambda}(h)\in Z(\H_n)$.
  \item The set $B=\{b_\alpha\mid\alpha\vdash n\}$ is a $\Q[\xi]$--basis for $Z(\H)$.\label{thm:part:norms.form.Qbasis}
  \end{enumerate}
\end{thm}

\begin{proof}
Part (1) is \cite[Lemma (2.12)]{Jon90}.
Part (2) is \cite[Proposition (2.13)]{Jon90}.
Part (3) is \cite[Theorem (3.33)]{Jon90}.
\end{proof}

Note that the set $B$ in Theorem \ref{thm:Jones.norm.properties}\eqref{thm:part:norms.form.Qbasis} is indexed by partitions of $n$ (not compositions).  See Proposition \ref{lem:norms.from.bips}.

\begin{thm}[\protect{\cite[(2.30)]{Jon90}}]\label{thm:norms.decomp.Jones.2.30}
Let $\lambda,\mu\vDash n$.  Let $S_\lambda$ and $S_\mu$ be parabolic subgroups of $S_n$.  If $b\in Z_\H(\H_\lambda)$ then
\[\N_{S_n,S_\lambda}\left(b\right)=\sum_{d\in\D_{\lambda\mu}}\N_{S_\mu,S_\lambda^d\cap S_\mu}\left(\T_{d^{-1}}b\T_d\right).\]
\end{thm}

{
\begin{prop}[\protect{\cite[(2.32)]{Jon90}}]\label{prop:norms.split.Jones.2.32}
Let $\lambda$ and $\mu$ be compositions satisfying $|\lambda|+|\mu|=n$, and let $\lambda'\le\lambda$, $\mu'\le\mu$.  {Let $x\in S_{\lambda'}$ and $y\in S_{\mu'}$.
Then
\[\N_{S_\lambda\times S_\mu,S_{\lambda'}\times S_{\mu'}}\left(\T_{xy}\right)=\N_{S_\lambda,S_{\lambda'}}\left(\T_{x}\right)
\N_{S_\mu,S_{\mu'}}\left(\T_{y}\right).\]}
\end{prop}
}

{
\begin{prop}[\protect{\cite[(3.29)]{Jon90}}]\label{prop:norms.specialized.Jones.3.29}
For $\alpha\vdash n$, we have
  \[b_\alpha|_{\xi=0}=\left[N_{S_n}(S_\alpha):S_\alpha\right]\underline{C_\alpha}.\]
\end{prop}
}

\subsection{The minimal basis}
The minimal basis is the analogue of the class sum basis for the {centre of the} group algebra.  Its existence was shown in \cite{GR97}, and it was explicitly described in \cite{Fmb}.

\begin{thm}[\cite{GR97}]\label{thm:geck.rouquier}
There exists a set of elements $\{\Gamma_\lambda\mid\lambda\vdash n\}\subseteq Z(\H)$ characterized by the properties
\begin{enumerate}
\item $\Gamma_\lambda|_{\xi=0}=\sum_{w\in C_\lambda}\T_w$, and
\item $\Gamma_\lambda-\sum_{w\in C_\lambda}\T_w$ contains no shortest elements of any conjugacy class.
\end{enumerate}
These elements form a $\Z[\xi]$--basis for $Z(\H)$.
\end{thm}

An element in $\H^+$ is said to be \emph{primitive} if, when written as a linear combination of $\{\T_w\mid w\in S_n\}$, its coefficients have no common factors over $\Z[\xi]$.  The main result of \cite{Fmb} is the following (using the partial order introduced in Section \ref{subsec:Hecke.algebras}):

\begin{thm}[\cite{Fmb}]The set $\{\Gamma_\lambda\mid\lambda\vdash n\}$ is the set of primitive minimal elements of $Z(\H)^+$.
\end{thm}

Also we have:
\begin{lem}\label{lem:class.elt.properties}\text{}
\begin{enumerate}
\item  $\Gamma_{(n)}=\eta_{(n)}$.
\item If $r\T_{w_\lambda}\le h\in Z(\H)^+$ then $r\Gamma_\lambda\le h$.
\end{enumerate}
\end{lem}
\begin{proof}
(1) is immediate from {Theorem} \ref{lem:etacox} (3) and the characterization of Theorem \ref{thm:geck.rouquier}.
{Part} (2) is \cite[Corollary 4.6]{Fmb}.
\end{proof}

\section{An inner product on the Hecke algebra}\label{sec:inner.prod}

The standard trace function $\tau$ on $\H$ is defined by $\tau(\T_w)=1$ if $w=1$ and $0$ otherwise.  A generalization of this trace was defined in \cite{Fcent} as follows.  Fix an element $h=\sum_{w\in S_n}r_w\T_w\in Z(\H)$, with $r_w\in\Z[\xi]$.  For $w\in S_n$, set $\h(\T_w):=r_w$, and extend linearly to all of $\H$.  If {$h=\T_1$} then this is simply $\tau$.  In this section we introduce a more flexible alternative formulation of this map as an inner product.

Define a map
$\l\ ,\ \r:\H\times\H\to\Z[\xi]$ by setting
\begin{equation}\label{eq:innerprod}
\l\sum_{w\in S_n}r_w\T_w,\sum_{w\in S_n}r'_w\T_w\r=\sum_{w\in S_n}r_wr'_w.
\end{equation}
If $h\in Z(\H)$, then  \eqref{eq:innerprod} defines the map $\h$, since $\l h,\T_w\r$ gives the coefficient of $\T_w$ in $h$.
This map and Proposition \ref{prop:innerprod} were developed during discussions between the first author and Leonard Scott in 1999.

\begin{prop}\label{prop:innerprod} Let $u,v,w\in S_n$.  The map defined in \eqref{eq:innerprod} is an inner product on $\H$ satisfying
\begin{equation}\l \T_u\T_v\T_{u^{-1}},\T_w\r=\l \T_v,\T_{u^{-1}}\T_w\T_u\r\label{eq:inner.prod.property}\end{equation}
\end{prop}
\begin{proof}
The properties of an inner product are easily verified from the definition.  We prove \eqref{eq:inner.prod.property} by induction on the length of $u$.

If $\ell(u)=1$, then we may set $u=s\in S$, and reduce the problem to considering $v$ and $w$ in the same $\l s\r$--$\l s\r$ double coset of $S_n$ (if they are in different double cosets then both sides of \eqref{eq:inner.prod.property} will be zero and the statement holds).  Further, if $v=w$ and $u=s$ then the statement holds by the symmetry of the inner product.  So to prove the result for $\ell(u)=1$ we need to check the cases where $v\neq w$ in the same $\l s\r$--$\l s\r$ double coset.

There are two cases for such a double coset $\l s\r d\l s\r$ in the symmetric group: either $ds=sd$ or $ds\neq sd$.  If $ds=sd$ then the double coset consists of only $\{d,ds\}$ and there is just one case to check:
\begin{align*}
\l \T_d,\T_s\T_{ds}\T_s\r &=\l \T_d,\T_s\left(\T_d+\xi\T_{ds}\right)\r=\l \T_d, { \T_{sd} } +\xi\T_d+\xi^2\T_{ds}\r\\
&=\xi=\l \T_d+\xi\T_{ds},\T_{ds}\r=\l\T_s\T_d\T_s,\T_{ds}\r.
\end{align*}
When $ds\neq sd$ the double coset has four elements $\{d,ds,sd,sds\}$, and so there are ${4\choose 2}=6$ cases:
\begin{align*}
\l \T_d,\T_s\T_{ds}\T_s\r &=\l \T_d,\T_{sd}+\xi\T_{sds}\r=0=\l \T_{sds},\T_{ds}\r=\l\T_s\T_d\T_s,\T_{ds}\r,\\
\l \T_d,\T_s\T_{sds}\T_s\r &=\l \T_d,\T_d+\xi\T_{ds}+\xi\T_{sd}+\xi^2\T_{sds}\r=1=\l \T_s\T_d\T_s,\T_{sds}\r,\\
\l \T_{ds},\T_s\T_{sd}\T_s\r &=\l \T_{ds},\T_{ds}+\xi\T_{sds}\r=1=\l\T_s\T_{ds}\T_s,\T_{sd}\r,\\
\l \T_{ds},\T_s\T_{sds}\T_s\r &=\xi=\l\T_s\T_{ds}\T_s,\T_{sds}\r.
\end{align*}
The cases
\[\l \T_d,\T_s\T_{sd}\T_s\r=\l\T_s\T_d\T_s,\T_{sd}\r\text{ and }\l \T_{sd},\T_s\T_{sds}\T_s\r =\l\T_s\T_{sd}\T_s,\T_{sds}\r\]
are symmetric to cases listed above.

Now let $u=u_1s$ for some $u_1\in S_n$ and $s\in S$ with $\ell(u_1s)=\ell(u_1)+1$.  We have
\begin{align*}
  \l \T_u\T_v\T_{u^{-1}},\T_{w}\r &=  \l \T_{u_1}\left(\T_{s}\T_v\T_{s}\right)\T_{u_1^{-1}},\T_{w}\r\\
  &= \l \T_{s}\T_v\T_{s},\T_{u_1^{-1}}\T_{w}\T_{u_1}\r\quad\text{(by linearity of the inner}\\
  &\hspace{2cm}\text{ product, and by induction)}\\
  &= \l \T_v,\T_{s}\left(\T_{u_1^{-1}}\T_{w}\T_{u_1}\right)\T_{s}\r\quad\text{(for the same reasons)}\\
  &= \l \T_v,\T_{u}^{-1}\T_{w}\T_{u}\r,
\end{align*}
which completes the proof.
\end{proof}

\begin{cor}\label{cor:innerprod.norm}
Let $\lambda,\mu\vDash n$.  Then
\[\l \N_{S_n,1}(\T_{w_\lambda}),\T_{w_\mu}\r =\l \N_{S_n,1}(\T_{w_\mu}),\T_{w_\lambda}\r .\]
\end{cor}
\begin{proof}
Immediate from \ref{prop:innerprod}, the symmetry of the inner product, and the fact that the norm from the identity sums over the entire group.
\end{proof}

The following Lemma was stated originally in terms of the function $\h$ for a fixed $h\in Z(\H)$.  The statement below is a direct translation of that result.
\begin{lem}[{\cite[Lemma {(3.5)}]{Fcent}}]\label{lem:innerprod:symm}
If $h\in Z(\H)$ then $\l h,\T_w\T_v\r=\l h,\T_v\T_w\r$.
\end{lem}

{A key property of the inner product in our context is the following:
\begin{lem}\label{lem:class.elt.coeff.in.central.elt}
Let $\lambda\vdash n$.  The coefficient of $\Gamma_\lambda$ in $h\in Z(\H)$ is $\l h,\T_{w_\lambda}\r$.
\end{lem}
\begin{proof}
  Immediate from Lemma \ref{lem:class.elt.properties} (2).
\end{proof}
}

\section{The norm of the identity in terms of the class elements}\label{sec:norm.of.one}

The goal of this section is to prove Theorem \ref{thm:norm.of.1}, which gives the coefficient of a class element in the norm of the identity.  This result is required for the proof of Theorem \ref{thm:coeff.cox.elt.in.balpha}.  Eventually, Theorem \ref{thm:norm.of.1} is subsumed in the statement of Theorem \ref{thm:coeff.Glambda.in.balpha}.

\begin{lem}\label{lem:innerprod:centralelt:godown}
For $1\le i\le n-1$,
\[\l \Gamma_{(n)},\T_{s_1s_2\dots s_ns_i\dots s_1}\r =\xi\l \Gamma_{(n)},\T_{s_1s_2\dots s_ns_{i-1}\dots s_1}\r .\]
\end{lem}
\begin{proof}
{Note that $s_1s_2\dots s_ns_{i-1}\dots s_1$ and $s_1s_2\dots s_ns_i\dots s_2$ are in the same conjugacy class and are conjugate by some $x\in S_n$ satisfying the Geck-Pfeiffer property II (Theorem 1.1(ii) of \cite{GP93}), and so the coefficients of their corresponding elements in any central element of the Hecke algebra are equal (using Lemma \ref{lem:innerprod:symm}).  In particular, $\l \Gamma_{(n)},\T_{s_1s_2\dots s_ns_{i-1}\dots s_1}\r=\l \Gamma_{(n)},\T_{s_1s_2\dots s_ns_i\dots s_2}\r$. }
Using {Lemma} \ref{lem:innerprod:symm} again, we have that
\begin{align*}
\l \Gamma_{(n)},\T_{s_1s_2\dots s_ns_i\dots s_1}\r
 &=\l \Gamma_{(n)},\T_{s_1}\T_{s_1s_2\dots s_ns_i\dots s_2}\r \\
 &=\l \Gamma_{(n)},\left(\T_1+\xi\T_{s_1}\right)\T_{s_2\dots s_ns_i\dots s_2}\r \\
 &=\l \Gamma_{(n)},\T_{s_2\dots s_ns_i\dots s_2}+\xi\T_{s_1s_2\dots s_ns_i\dots s_2}\r \\
 &=\l \Gamma_{(n)},\xi\T_{s_1s_2\dots s_ns_i\dots s_2}\r \text{ (by linearity, }\\
 &\hspace{2cm}\text{Theorem \ref{lem:etacox}(2) and Lemma \ref{lem:class.elt.properties}(1))} \\
 &=\l \Gamma_{(n)},\xi\T_{s_1s_2\dots s_ns_{i-1}\dots s_1}\r ,
\end{align*}
{which proves the lemma.}
\end{proof}

\begin{lem}\label{lem:mountain}
For $n\ge 2$, the coefficient of $\T_{s_1s_2\dots s_n\dots s_2s_1}$ in $\Gamma_{(n)}$ is $\xi^{n-2}$.
\end{lem}
\begin{proof}
This follows from repeated application of Lemma \ref{lem:innerprod:centralelt:godown}.
\end{proof}

\begin{thm}\label{thm:norm.of.1}
The coefficient of $\Gamma_\lambda$ in $\N_{S_n,1}(\T_1)$ is $[S_n:S_\lambda]\xi^{l_\lambda}$.  Equivalently,
\[\l \N_{S_n,1}(\T_1),\T_{w_\lambda}\r =[S_n:S_\lambda]\xi^{l_\lambda}\]
where $w_\lambda$ is as defined in Section \ref{subsec:symm.group}.
\end{thm}
\begin{proof}
Firstly note that the two statements in the theorem are in fact equivalent by {Lemma \ref{lem:class.elt.coeff.in.central.elt}}.

Secondly note that the symmetry of the norm from the trivial subgroup with respect to the inner product (Corollary \ref{cor:innerprod.norm}) means that
\[\l \N_{S_n,1}(\T_1),\T_{w_\lambda}\r =\l \N_{S_n,1}(\T_{w_\lambda}),\T_1\r { . } \]
{Thus} the problem reduces to identifying the coefficient of $\T_1$ in the norm of {$\T_{w_\lambda}$} from the trivial subgroup to $S_n$.

The transitivity of the norm implies
\begin{align}
\N_{S_n,1}(\T_{w_\lambda})
&=\N_{S_n,S_\lambda}(\N_{S_\lambda,S_{\lambda-1}}(\N_{S_{\lambda-1},1}(\T_{w_\lambda})))\notag\\
&=\N_{S_n,S_\lambda}(\N_{S_\lambda,S_{\lambda-1}}(\eta_\lambda))\notag\\
&=\N_{S_n,S_\lambda}(\eta_\lambda \N_{S_\lambda,S_{\lambda-1}}(\T_1)).\label{eq:norm.of.one}
\end{align}
The result therefore depends on finding the coefficient of $\T_1$ in expression \eqref{eq:norm.of.one}.  Now \eqref{eq:norm.of.one} is a linear combination of terms of {the} form $\T_{d^{-1}}\T_w\T_d$ where $d$ is a distinguished right coset representative of $S_\lambda$ in $S_n$, and {where} $\T_w$ occurs in $\eta_\lambda \N_{S_\lambda,S_{\lambda-1}}(\T_1)$ {(and is therefore an element of $\H_\lambda$)}.
Since $d$ is a coset representative, $\T_1$ occurs in $\T_{d^{-1}}\T_w\T_d$ only if $w=1$ (by Theorem \ref{thm.Shi.Bruhat}).

Similarly, by Theorem \ref{thm.Shi.Bruhat}, it is straightforward to see that $\T_1$ occurs in $\eta_\lambda \N_{S_\lambda,S_{\lambda-1}}(\T_1)$ only when $\T_{w^{-1}}$ occurs in $\N_{S_\lambda,S_{\lambda-1}}(\T_1)$ for some $\T_w$ occurring in $\eta_\lambda$.

The norm $\N_{S_\lambda,S_{\lambda-1}}(\T_1)$ splits into commuting factors along the components of $\lambda$.  That is, if $\lambda=(\lambda_1,\dots,\lambda_r)$ then
\begin{equation}
\N_{S_\lambda,S_{\lambda-1}}(\T_1)=\prod_{i=1}^r\N_{S_{\lambda_i},S_{\lambda_i-1}}(\T_1).\label{eq:norm.lambda.of.1.factors}
\end{equation}
The non-trivial {right} coset representatives of each $S_{\lambda_i-1}$ in $S_{\lambda_i}$ are of {the} form {$s_{\lambda_1+\dots+\lambda_{i-1}-1}\dots s_{\lambda_1+\dots+\lambda_{i-1}+\lambda_i+j}$} for $1\le j\le \lambda_i-1$.  Thus each factor in \eqref{eq:norm.lambda.of.1.factors} is a sum of products of {the} form
\begin{equation}\label{eq:slant.product}
\T_{s_{x+j}\dots s_{x+\lambda_i-1}}\T_{s_{x+\lambda_i-1}\dots s_{x+j}}
\end{equation}
where $x=\lambda_1+\dots+\lambda_{i-1}$.
{By Lemma \ref{lem:slant.product}, each non-identity term in the expansion of \eqref{eq:slant.product} is of the form $\xi\T_w$ where $w$ is a transposition.  }
Thus the non-identity terms in $\N_{S_\lambda,S_{\lambda-1}}(\T_1)$ are products of terms of the form $\xi\T_w$ where $w$ is a transposition in $S_\lambda$.
In addition, if $w$ is the longest transposition in $S_\lambda$, then $\xi^t$ is the coefficient of $\T_w$ where $t=|\{\lambda_i\ge 2\}|$, by \eqref{eq:norm.lambda.of.1.factors}.

Therefore, $\eta_\lambda \N_{S_\lambda,S_{\lambda-1}}(\T_1)$ {contains $\T_1$ only when $\eta_\lambda$ contains a transposition}.
{By Lemma \ref{lem:etacox} {(2)}, $\eta_\lambda$ contains only the longest transposition in $S_\lambda$, and its coefficient by Lemma \ref{lem:mountain} is} $\prod_{i=1}^t\xi^{\lambda_i-2}=\xi^{\lambda_1+\dots+\lambda_t-2t}$.
Hence, the coefficient of $\T_1$ in $\eta_\lambda \N_{S_\lambda,S_{\lambda-1}}(\T_1)$ is $\xi^{\lambda_1+\dots+\lambda_t-2t}\xi^t=\xi^{l_\lambda}$.

So $\T_1$ occurs in $\N_{S_n,1}(\T_{w_\lambda})$ for each coset representative of $S_\lambda$ in $S_n$.  The statement follows.
\end{proof}

\section{A partial order on the norm basis and the coefficient of the Coxeter class}\label{sec:partial.order.on.norms}

In this section we prove that the elements of the norm basis {ordered by refinement of partitions are also} ordered by the Hecke algebra order defined in Section \ref{subsec:Hecke.algebras}.
While this result might be of independent interest,
our chief use for this fact in the present paper is to obtain explicitly the coefficient of the Coxeter class element in a given norm basis element (see Theorem \ref{thm:coeff.cox.elt.in.balpha}).

Recall that the set $B:=\{b_\alpha\mid \alpha\vdash n\}$ is the $\Q[\xi]$--basis for the centre $Z(\H)$ defined in Section \ref{subsec:norm.basis}.

\begin{lem}\label{lem:order.on.balpha.max.sgp}
For any {integers} $n$ and $k$ with $n-k\le k<n$, {$\xi b_{(n)}< b_{(k,n-k)}$}.
\end{lem}
\begin{proof}
Using the fact that the double coset representative $d_{n-k}$ as defined in Section \ref{sec.double.cosets} is also a right coset representative of $S_{(k,n-k)}$ in $S_n$, we have
\begin{align*}
b_{(k,n-k)}&=\N_{S_n,S_{(k,n-k)}}(\eta_{(k,n-k)})\\
           &\ge\T_{d_{n-k}}\T_{s_1\dots s_{k-1}s_{k+1}\dots s_{n-1}}\T_{d_{n-k}}\\
           &=\T_{s_1\dots s_{k-1}s_{k+1}\dots s_{n-1}}\T_{d_{n-k}}^2\quad\text{(since $d_{n-k}$ commutes }\\
           &\hspace{2cm}\text{with $s_1\dots s_{k-1}s_{k+1}\dots s_{n-1}$)}\\
           &>\T_{s_1\dots s_{k-1}s_{k+1}\dots s_{n-1}}(\xi\T_{s_k})\quad\text{(since {$\xi\T_{s_k}\le \T_{d_{n-k}}^2$} )}\\
           &=\xi\T_{s_1\dots s_{k-1}s_{k+1}\dots s_{n-1}}\T_{s_k}.
\end{align*}
Since $s_1\dots s_{k-1}s_{k+1}\dots s_{n-1}s_k$ is a minimal element of the Coxeter class $C_{(n)}$, and since $b_{(k,n-k)}\in Z(\H)^+$ (\cite[Proposition (5.3)(ii)]{Fmb}), it follows that $\xi\Gamma_{(n)}\le b_{(k,n-k)}$ by \cite[Corollary (4.6)]{Fmb}.  Since $b_{(n)}=\Gamma_{(n)}$, and since $\T_{w_{(k,n-k)}}$ occurs in $b_{(k,n-k)}$ but not in $b_{(n)}$, the Lemma follows.
\end{proof}

\begin{thm}\label{thm:balpha.partial.order}
Let $\lambda,\mu\vDash n$.  If $\lambda<\mu$ in the refinement order, then $\xi^{l_\mu-l_\lambda}b_\mu< b_\lambda$.
\end{thm}
\begin{proof}
{We first} prove the {theorem when} $\lambda=(\lambda_1,\dots,\lambda_r)$ and $\mu=(\lambda_1,\dots,\lambda_{r-2},\lambda_{r-1}+\lambda_r)$.
{In this case we have}
\begin{align*}
b_\lambda
&=\N_{S_n,S_{\lambda}}(\eta_{\lambda})\\
&=\N_{S_n,S_{\mu}}\left(\N_{S_{\mu},S_{\lambda}}(\eta_{\lambda})\right)\\
&=\N_{S_n,S_{\mu}}\left(\left(\prod_{i=1}^{r-2}\N_{S_{\lambda_i},S_{\lambda_i}}(\eta_{\lambda_i})\right) \N_{S_{\lambda_{r-1}+\lambda_r},S_{(\lambda_{r-1},\lambda_r)}}(\eta_{(\lambda_{r-1},\lambda_r)})\right)\\
&>\N_{S_n,S_{\mu}}\left(\left(\prod_{i=1}^{r-2}\N_{S_{\lambda_i},S_{\lambda_i}}(\eta_{\lambda_i})\right)
\xi\N_{S_{\lambda_{r-1}+\lambda_r},S_{\lambda_{r-1}+\lambda_r}}(\eta_{(\lambda_{r-1}+\lambda_r)})\right)\\
&\hspace{3cm}\text{(by Lemma \ref{lem:order.on.balpha.max.sgp})}\\
&=\xi\N_{S_n,S_{\mu}}\left(\N_{S_{\mu},S_{\mu}}(\eta_{\mu})\right)\\
&=\xi\N_{S_n,S_{\mu}}\left(\eta_{\mu}\right)\\
&=\xi b_{\mu}.
\end{align*}

{Similar arguments demonstrate the validity of the theorem when $\mu=(\lambda_1,\lambda_2,\dots , \lambda_{i-1}, \lambda_{i}+\lambda_{i+1}, \lambda_{i+2}, \dots , \lambda_r)$ with $1\le i \le r-2.$} The extension to arbitrary $\lambda<\mu$ is {then} immediate.
\end{proof}

\begin{cor}\label{cor:balpha.class.elts.partial.order}
Write $b_\alpha=\sum_{\lambda\vdash n}{r_{\alpha,\lambda}}\Gamma_\lambda$ for $\alpha\vdash n$ and $r_{\alpha,\lambda}\in\Z[\xi]$.  Then $\alpha<\beta$ implies $\xi r_{\beta,\lambda}\le r_{\alpha,\lambda}$ for any $\lambda\vdash n$
\end{cor}
\begin{proof}
The partial order of Theorem \ref{thm:balpha.partial.order} implies that for an arbitrary basis element $\T_w$ of $\H$, if $\lambda<\mu$ then the coefficients $r_{w,\lambda}$ and $r_{w,\mu}$ of $\T_w$ in $b_\lambda$ and $b_\mu$ respectively satisfy the relation $\xi r_{w,\mu}\le r_{w,\lambda}$.  In particular, this is the case when $w$ is a shortest element of a conjugacy class.  It follows that for any sequence of partitions totally ordered by refinement  between $(1^n)$ and $(n)$, the coefficients of any given class element also satisfy the inequality {by Lemma \ref{lem:class.elt.coeff.in.central.elt}}.
\end{proof}

\begin{thm}\label{thm:coeff.cox.elt.in.balpha}
The coefficient of the Coxeter class element $\Gamma_{(n)}$ in $b_\alpha$ is $\xi^{n-1-l_\alpha}$.
\end{thm}
\begin{proof}
We have from Lemma \ref{lem:class.elt.properties}{(1)} that the coefficient of $\Gamma_{(n)}$ in $b_{(n)}$ is 1.  We also have that the coefficient of $\Gamma_{(n)}$ in $b_{(1^n)}$ is $\xi^{n-1}$ by Theorem \ref{thm:norm.of.1}.  Clearly both {of} these coefficients satisfy the theorem statement.  Given any sequence of partitions $\alpha_i$ totally ordered by refinement between $(1^n)$ and $(n)$:
\[(1^n)=\alpha_0<\dots<\alpha_{n-1}=(n),\]
Corollary \ref{cor:balpha.class.elts.partial.order} gives the following relations, where $r_{i}$ is the coefficient of $\Gamma_{(n)}$ in $b_{\alpha_i}$:
\[\xi^{n-1}=r_{0}\ge\xi r_{1}\ge\dots\ge \xi^{i}r_i\ge\dots\ge\xi^{n-1}r_{n-1}=\xi^{n-1}.
\]
It follows that all the terms in the above order are equal, and $\xi^ir_i=\xi^{n-1}$, so that $r_i=\xi^{n-1-i}$.  Since length increases by one with each step of the refinement order, $i=l_{\alpha_i}$. The result follows since every partition is part of such a full sequence of partitions between $(1^n)$ and $(n)$ ordered by refinement with length increasing by one.
\end{proof}
\begin{rem}
Applying the map from $\H\to \H_q$ (as described in the Remark in Section \ref{subsec:Hecke.algebras}) to the expansion of $b_\alpha$, the coefficient of the Coxeter class element in Theorem \ref{thm:coeff.cox.elt.in.balpha} in terms of the $\Z[q,q^{-1}]$--basis for $\H_q$ is $\left(q^{-1}(q-1)\right)^{n-1-l_\alpha}$.
Note that under the map $\H\to \H_q$, $\Gamma_{(n)}\mapsto q^{-\frac{n-1}{2}}\Gamma_{(n),q}$ and $b_\alpha\mapsto q^{-\frac{l_\alpha}{2}}b_{\alpha,q}$, where $\Gamma_{(n),q}$ and $b_{\alpha,q}$ are the corresponding $\H_q$ versions.
This is an example of a result whose statement over $\Z[q,q^{-1}]$ is independent of the
partial order on $\H^+$ --- which makes sense only over $\Z[\xi]$ --- but whose proof is
made possible by treating the Hecke algebra as a $\Z[\xi]$--module and using the partial
order.
\end{rem}

\section{Projections of norms onto maximal parabolic subalgebras}\label{sec:proj.norms.onto.maximal}

The main result of this section is a projection formula for $b_\alpha$ (a norm basis element) onto a maximal parabolic subgroup.  Results {in} this section are from \cite{Jon87} but have not appeared in the literature.

Let $B=\{b_\alpha\mid\alpha\vdash n\}$ (see Section \ref{sec:bases.for.centre}).
If $\theta=(\theta_1,\dots,\theta_t)$ is a multipartition of $n$, write $\eta_\theta:=\eta_{\theta_1}\dots\eta_{\theta_t}$.

For $\lambda\vDash n$, define the projection $\pi_\lambda:\H_n\to\H_\lambda$ by setting
\[\pi_\lambda(\T_w)=\begin{cases}\T_w&\text{ if }w\in S_\lambda\\ 0&\text{ otherwise.}\end{cases}\]
and extending to $\H_n$ linearly.

\begin{thm}[\cite{Jon87}]\label{thm:proj.balpha.to.max.sgp}
Let $b_{\alpha}\in B$ and write $\pi:=\pi_{(k,n-k)}$.  Then
\begin{equation}\label{eq:proj.balpha.to.max.sgp}
\pi (b_{\alpha})=\sum_{(\mu,\nu)\in\Lambda_{(k,n-k)}(\alpha)}z_{\mu,\nu}\N_{S_{(k,n-k)},S_{(\mu,\nu)}}(\eta_{(\mu,\nu)}),
\end{equation}
with the $z_{\mu,\nu}$ nonnegative integers.
\end{thm}

\begin{proof} The proof is by induction on $n.$ The theorem is easily verified for $n=2$, so assume it in all cases less than $n$.

For $(\mu,\nu)\in\Lambda_{(k,n-k)}(\alpha)$ we have
\begin{align*}
b_{\alpha}&=\N_{S_{n},S_{\alpha}}(\eta_{\alpha})\\
&=\N_{S_{n},S_{(k,n-k)}}\left(\N_{S_{(k,n-k)},S_{(\mu,\nu)}}(\eta_{(\mu,\nu)})\right)\\
&\hspace{4cm}\text{(by Theorem \ref{thm:Jones.norm.properties} and Proposition \ref{lem:norms.from.bips})}\\
&=\sum_{m=0}^{n-k}\N_{S_{(k,n-k)},P_m}\left({\T_{d_{m}}\N_{S_{(k,n-k)},S_{(\mu,\nu)}}(\eta_{(\mu,\nu)})\T_{d_{m}}}\right)\\
&\hspace{4cm}\text{(by {Theorem \ref{thm:norms.decomp.Jones.2.30}}) } \\
&=\sum_{m=0}^{n-k}\N_{S_{(k,n-k)},P_m}\left({\T_{d_{m}}\N_{S_k,S_\mu}(\eta_\mu)
\N_{S_{n-k},S_\nu}(\eta_\nu)\T_{d_{m}}}\right )\\
&\hspace{4cm}\text{(by {Proposition \ref{prop:norms.split.Jones.2.32}.}) }
\end{align*}

Consider the general term in this last sum. For this term, write $\N_{S_k,S_\mu}(\eta_{\mu})$ as $A_{k}+B_{k}$ where $A_{k}$ is the projection of $\N_{S_{k},S_\mu}(\eta_\mu)$ onto $\H_{(k-m,m)}$ in $\H_k$.  Similarly, write $\N_{S_{n-k},S_\nu}(\eta_\nu)$ as $A_{n-k}+B_{n-k}$ where $A_{n-k}$ is the projection of $\N_{S_{n-k},S_\nu}(\eta_\nu)$ onto $\H_{(m,n-k-m)}$ in $\H_{n-k}$.
Thus
\[\N_{S_k,S_\mu}(\eta_\mu)\N_{S_{n-k},S_\nu}(\eta_\nu)=A_kA_{n-k}+A_kB_{n-k}+B_kA_{n-k}+B_kB_{n-k}.\]
By induction,
\[A_k=\sum_{(\mu_1,\mu_2)\in\Lambda_{(k-m,m)}(\mu)}z_{\mu_1,\mu_2}\N_{S_{(k-m,m)},S_{(\mu_1,\mu_2)}}(\eta_{(\mu_1,\mu_2)}),\]
and
\[A_{n-k}=\hspace{-2mm}\sum_{(\nu_1,\nu_2)\in\Lambda_{(m,n-k-m)}(\nu)}\hspace{-2mm}z_{\nu_1,\nu_2}\N_{S_{(m,n-k-m)},S_{(\nu_1,\nu_2)}}(\eta_{(\nu_1,\nu_2)}),\]
with the $z_{\mu_1,\mu_2}$ and $z_{\nu_1,\nu_2}$ non-negative integers.
{We abbreviate the $(k-m,m,m,n-k)$--multipartition} $(\mu_1,\mu_2,\nu_1,\nu_2)$ of $\alpha$ by $\bar\alpha$.
Then
\begin{equation}\label{eq:Ak.An-k}
A_{k}A_{n-k}=\hspace{-5mm}\sum_{\substack{(\mu_1,\mu_2)\in\Lambda_{(k-m,m)}(\mu)\\ (\nu_1,\nu_2)\in\Lambda_{(m,n-k-m)}(\nu)}}\hspace{-5mm}z_{\mu_1,\mu_2}z_{\nu_1,\nu_2}
\N_{P_m,S_{\bar\alpha}}(\eta_{\bar\alpha})
\end{equation}
by splicing the norms back together utilizing {Proposition \ref{prop:norms.split.Jones.2.32}}.

Each of the terms in \eqref{eq:Ak.An-k} is then multiplied by $\T_{d_m}$ on the left and right. Write
\[
\N_{P_m,S_{\bar\alpha}}(\eta_{\bar\alpha})
=\sum_{r\in R}\T_{r^{-1}}\eta_{\bar\alpha}\T_r
\]
where $R$ is the set of distinguished right coset representatives for $S_{\bar\alpha}$ in $P_{m}$.

Because $d_m$ is a distinguished double coset representative of $S_{(k,n-k)}$ in $S_n$, we have $\ell(rd_m)=\ell(d_mr^{-1})=\ell(d_m)+\ell(r)$ for any $r\in R$.  Then, since $\hat{r}=d_mrd_m$ is a distinguished right coset representative for $S_{\bar\alpha}^{d_m}$ in $P_m$ (Corollary \ref{cor:drd.is.coset.rep}), {it follows that:}
\begin{align*}
\T_{d_m}\T_{r^{-1}}\eta_{\bar\alpha}\T_r\T_{d_m}
    &=\T_{d_mr^{-1}}\eta_{\bar\alpha}\T_{rd_m}\\
    &=\T_{\hat{r}^{-1}d_m}\eta_{\bar\alpha}\T_{d_m\hat{r}}\\
    &=\T_{\hat{r}^{-1}}\T_{d_{m}}\eta_{\bar\alpha}\T_{d_m}\T_{\hat{r}}.
\end{align*}
{Since $\eta_{\bar\alpha}\in \H_{\bar\alpha}$ and $S_{\bar\alpha}\le P_m$}, Corollary \ref{cor:tau.in.Tdm.Tgam.Tdm.3.17} implies that
\[
\T_{d_m}\eta_{\bar\alpha}\T_{d_m}=\eta_{\bar\alpha^{d_m}}+\sum_{w\in S_n}f_{w}\T_w,\quad\text{where $\# (w)\ge 1$ if $f_{w}\ne 0$.  }
\]
Since $\# (\hat{r})=0$, {we have from} Corollary \ref{lem:tau.in.prod.with.Bruhat} {(2) that}
\begin{multline*}
\T_{\hat{r}^{-1}}\T_{d_m}\eta_{\bar\alpha}\T_{d_m}\T_{\hat{r}}\\
=\T_{\hat{r}^{-1}}\eta_{\bar\alpha^{d_m}}\T_{\hat{r}}+\sum_{w\in S_n}g_{w}\T_w,\quad\text{where $\# (w)\ge 1$ if $g_{w}\ne 0$,}
\end{multline*}
and so,
\begin{multline*}
\T_{d_m}\N_{P_m,S_{\bar\alpha}}(\eta_{\bar\alpha})\T_{d_m}\\
=\N_{P_m,S_{\bar\alpha^{d_m}}}(\eta_{\bar\alpha^{d_m}})+\sum_{w\in S_n}h_{w}\T_w,
\quad\text{where $\# (w)\ge 1$ if {$h_{w}\ne 0$}.}
\end{multline*}
Thus, applying Theorem \ref{thm:Jones.norm.properties} \eqref{thm:part:norm.transitivity} to $\H_{(k,n-k)}$,
\begin{multline}\label{eq:proj.norms.AkAn-k}
\pi\left(\N_{S_{(k,n-k)},P_m}\left(\T_{d_m}A_kA_{n-k}\T_{d_m}\right)\right)\\
=\sum_{\substack{(\mu_1,\mu_2)\in\Lambda_{(k-m,m)}(\mu)\\ (\nu_1,\nu_2)\in\Lambda_{(m,n-k-m)}(\nu)}}z_{\mu_1,\mu_2}z_{\nu_1,\nu_2}
\N_{S_{(k,n-k)},S_{\bar\alpha}^{d_m}}(\eta_{\bar\alpha^{d_m}}).
\end{multline}
Collecting common terms and re-indexing, \eqref{eq:proj.norms.AkAn-k} may be written as
\begin{multline*}
\pi\left(\N_{S_{(k,n-k)},P_m}\left(\T_{d_m}A_kA_{n-k}\T_{d_m}\right)\right)\\
=\sum_{(\mu,\nu)\in\Lambda_{(k,n-k)}(\alpha)}z_{\mu,\nu}\N_{S_{(k,n-k)},S_{(\mu,\nu)}}(\eta_{(\mu,\nu)}),
\end{multline*}
with the $z_{\mu,\nu}$ non-negative integers.

We now examine the product $\T_{d_{m}}B_{k}B_{n-k}\T_{d_{m}}.$ Since $B_{k}$ projects to zero on $\H_{(k-m,m)}$, every nonzero term of $B_{k}$ contains $s_{k-m}$.  Similarly, every nonzero term of $B_{n-k}$ contains $s_{k+m}$.  Let $\T_{w_{k}}$ be a nonzero term in $B_{k}$ and $\T_{w_{n-k}}$ a nonzero term in $B_{n-k}$, assuming for simplicity that the coefficients are one.
Now $w_{k}w_{n-k}\in S_{(k,n-k)}$ with $\ell (w_{k}w_{n-k})=\ell (w_{k})+\ell (w_{n-k})$, and since $d_m$ is a distinguished double coset representative of $S_{(k,n-k)}$ in $S_n$ we have
\begin{equation}\label{eq:Tdm.on.Bi}
\T_{d_{m}}\T_{w_{k}}\T_{w_{n-k}}\T_{d_{m}}=\T_{d_{m}}\T_{w_{k}w_{n-k}d_{m}}.
\end{equation}
If $d_{m}w_{k}w_{n-k}d_{m}\in S_{(k,n-k)}$, then $d_{m}w_{k}w_{n-k}d_{m}\in P_{m}$ {(by definition of $P_m$),} which implies that $w_{k}w_{n-k}\in P_{m}^{d_{m}}=P_{m}.$ This is a contradiction since $w_{k}$ contains $s_{k-m}\not\in P_{m}$ or since $w_{n-k}$ contains $s_{k+m}\not\in P_{m}.$ Thus, Corollary \ref{lem:coeff.tau.in.H.with.length.3.15}{(1)} implies that every term in \eqref{eq:Tdm.on.Bi} contains $s_k$ since $\# (d_{m}w_{k}w_{n-k}d_{m})>0$.  Hence, every term of $\T_{d_{m}}B_{k}B_{n-k}\T_{d_{m}}$ contains $s_k$.

Similarly, every term of the products $\T_{d_{m}}B_{k}A_{n-k}\T_{d_{m}}$ and $\T_{d_{m}}A_{k}B_{n-k}\T_{d_{m}}$ contains $s_k$.  Applying Corollary \ref{lem:tau.in.prod.with.Bruhat}{(2),} we conclude that
\[
\pi\left(\N_{S_{(k,n-k)},P_{m}}\left({\T_{d_m}(B_kA_{n-k}+A_kB_{n-k}+B_kB_{n-k})\T_{d_m}}\right)\right)=0.
\]
This completes the proof of the theorem.
\end{proof}

\begin{cor}\label{cor:coeff.norms.in.proj.to.max.sgp}
For each $(\mu,\nu)\in\Lambda_{(k,n-k)}$,
\[
z_{\mu,\nu}=[N_{S_n}(S_{(\mu,\nu)}):N_{S_{(k,n-k)}}(S_{(\mu,\nu)})].
\]
\end{cor}
\begin{proof}
We have from Theorem \ref{thm:proj.balpha.to.max.sgp} that
\[\pi (b_{\alpha})=\sum_{(\mu,\nu)\in\Lambda_{(k,n-k)}(\alpha)}z_{\mu,\nu}\N_{S_{(k,n-k)},S_{(\mu,\nu)}}(\eta_{(\mu,\nu)}).\]
By {Proposition \ref{prop:norms.specialized.Jones.3.29}},
when we specialize $b_\alpha$ by setting {$\xi=0$} we obtain $[N_{S_n}(S_\alpha):S_\alpha]\underline{C_\alpha}$.  When we project $\underline C_\alpha$ via $\pi$ we get a sum of disjoint subsets of $C_\alpha$, corresponding to $(k,n-k)$--multipartitions $(\mu,\nu)$ of $\alpha$.  Thus
\[\pi(b_{\alpha})|_{\xi=0}=\sum_{(\mu,\nu)\in\Lambda_{(k,n-k)}(\alpha)}[N_{S_n}(S_{(\mu,\nu)}):S_{(\mu,\nu)}]\underline{C_{(\mu,\nu)}}.\]
On the other hand
{\[
\N_{S_{(k,n-k)},S_{(\mu,\nu)}}(\eta_{(\mu,\nu)})|_{\xi=0}=[N_{S_{(k,n-k)}}(S_{(\mu,\nu)}):S_{(\mu,\nu)}]\underline{C_{(\mu,\nu)}}.
\]}
{Consequently,} for each $(\mu,\nu)\in\Lambda_{(k,n-k)}(\alpha)$ we have
\[
[N_{S_n}(S_{(\mu,\nu)}):S_{(\mu,\nu)}]=z_{\mu,\nu}[N_{S_{(k,n-k)}}(S_{(\mu,\nu)}):S_{(\mu,\nu)}],
\]
and the Corollary follows.
\end{proof}

\section{Coefficients of class elements in norms}\label{sec:coeffs.of.class.elts}

In this section we generalize Theorem \ref{thm:proj.balpha.to.max.sgp} to arbitrary parabolic subalgebras, and use this generalization to prove Theorem \ref{thm:coeff.Glambda.in.balpha} --- our main result --- on the coefficients of class elements in norms.

\begin{thm}\label{thm:proj.balpha.trickle}
Let $\lambda\vDash n$ and $\alpha\vdash n$.  Then
  \begin{equation}\label{eq:trickle.down.norms}
  \pi_\lambda(b_\alpha)=\frac{|S_n|}{|S_\lambda|\cdot|C_\alpha|_{S_n}}
\sum_{\theta\in\Lambda_\lambda(\alpha)}|C_\theta|_{S_\lambda}\N_{S_\lambda,S_\theta}(\eta_\theta).
  \end{equation}
\end{thm}

\begin{proof}
If $\lambda$ has one component, $\lambda=(n)$, then there is only one multipartition $\theta$ of shape $\lambda$ using the components {of} $\alpha$, namely $\theta=(\alpha)$ itself, and the statement is trivial.  If $\lambda=(k,n-k)$ then we can apply Theorem \ref{thm:proj.balpha.to.max.sgp} and Corollary \ref{cor:coeff.norms.in.proj.to.max.sgp} to obtain
\[\pi_\lambda(b_\alpha)=\sum_{\theta\in\Lambda_\lambda(\alpha)}\left[N_{S_n}(S_{\theta}):N_{S_\lambda}(S_\theta)\right]\N_{S_\lambda,S_\theta}(\eta_{\theta}).\]
An elementary calculation gives that
\[\left[N_{S_n}(S_{\theta}):N_{S_\lambda}(S_\theta)\right]=\frac{|S_n|\cdot|C_\theta|_{S_\lambda}}{|S_\lambda|\cdot|C_\alpha|_{S_n}},\]
and the statement for the case in which $\lambda$ has two components follows.

Now suppose inductively for some fixed $\lambda=(\lambda_1,\dots,\lambda_r)\vDash n$, with $r>2${,}
that \eqref{eq:trickle.down.norms} holds for all compositions of $n$ with fewer components than $\lambda$. In
particular, we assume that \eqref{eq:trickle.down.norms} holds for
$\lambda'=(\lambda_1,\dots,\lambda_{r-2},\lambda_{r-1}')$ where $\lambda_{r-1}'=\lambda_{r-1}+\lambda_r$.

Write $\theta'=\left(\theta_1,\dots,\theta_{r-2},\theta'_{r-1}\right)\in\Lambda_{\lambda'}(\alpha)$ with $\theta_i\vdash\lambda_i$ for $1\le i\le r-2$ and $\theta'_{r-1}\vdash\lambda'_{r-1}$. By our induction hypothesis we have
\begin{align}\label{eq:proj.b.lambda'}
&\pi_{\lambda'}(b_\alpha)\notag\\
&\ =\frac{|S_n|}{|S_{\lambda'}|\cdot|C_\alpha|_{S_n}}\sum_{\theta'\in\Lambda_{\lambda'}(\alpha)}|C_{\theta'}|_{S_{\lambda'}}\N_{S_{\lambda'},S_{\theta'}}(\eta_{\theta'})\notag\\
\begin{split}
&\ =\frac{|S_n|}{|S_{(\lambda_1,\dots,\lambda_{r-2})}|\cdot|S_{\lambda'_{r-1}}|\cdot|C_\alpha|_{S_n}}\\
&\quad\times\sum_{\theta'\in\Lambda_{\lambda'}(\alpha)}\left(|C_{(\theta_1,\dots,\theta_{r-2})}|_{S_{(\lambda_1,\dots,\lambda_{r-2})}}\cdot |C_{\theta'_{r-1}}|_{S_{\lambda'_{r-1}}}\right.\\
&\quad\left.\cdot\N_{S_{(\lambda_1,\dots,\lambda_{r-2})},S_{(\theta_1,\dots,\theta_{r-2})}}\left(\eta_{(\theta_1,\dots,\theta_{r-2})}\right) \cdot\N_{S_{\lambda'_{r-1}},S_{\theta'_{r-1}}}\big(\eta_{\theta'_{r-1}}\big)\right).
\end{split}
\end{align}

Let $\pi_{\lambda''}$ be the projection from the algebra $\H_{\lambda'}$ onto $\H_\lambda$.  Then
$\pi_{\lambda''}$ acts as the identity on $\H_{(\lambda_1,\dots,\lambda_{r-2})}$, and acts on norms
in $\H_{\lambda'_{r-1}}$ according to equation \eqref{eq:proj.balpha.to.max.sgp} and by the inductive hypothesis according to \eqref{eq:trickle.down.norms}.

Then $ \pi_\lambda=\pi_{\lambda''}\pi_{\lambda'}$, and $\pi_{\lambda''}$ acts
trivially on all but the last factor in each summand of \eqref{eq:proj.b.lambda'}, namely
$\N_{S_{\lambda'_{r-1}},S_{\theta'_{r-1}}}(\eta_{\theta'_{r-1}})$.
The image of this norm under the projection $\pi_{\lambda''}$ to the maximal parabolic subgroup
$S_{(\lambda_{r-1},\lambda_r)}$ of $S_{\lambda_{r-1}+\lambda_r}$ can thus be decomposed {by induction as} follows:
\begin{align*}
&\pi_{\lambda''}(\N_{S_{\lambda'_{r-1}},S_{\theta'_{r-1}}}(\eta_{\theta'_{r-1}}))\\
\begin{split}&\quad=\frac{|S_{\lambda'_{r-1}}|}{|S_{(\lambda_{r-1},\lambda_r)}|\cdot|C_{\theta'_{r-1}}|_{S_{\lambda'_{r-1}}}}\\
             &\qquad\times\sum_{\substack{(\theta_{r-1},\theta_r)\in\\ \Lambda_{(\lambda_{r-1},\lambda_r)}(\theta'_{r-1})}}
 |C_{(\theta_{r-1},\theta_r)}|_{S_{(\lambda_{r-1},\lambda_r)}} \cdot\N_{S_{(\lambda_{r-1},\lambda_r)},S_{(\theta_{r-1},\theta_r)}}(\eta_{(\theta_{r-1},\theta_r)}).\end{split}
\end{align*}
Composing $\pi_{\lambda''}$ and $\pi_{\lambda'}$ to give $\pi_\lambda$ we have (after some initial cancellations):
\begin{align*}
&\pi_\lambda(b_\alpha)\\
&=\frac{|S_n|}{|S_{(\lambda_1,\dots,\lambda_{r-2})}|\cdot|C_\alpha|_{S_n}}\\
&\quad\times\sum_{\theta'\in\Lambda_{\lambda'}(\alpha)}
    \frac{|C_{(\theta_1,\dots,\theta_{r-2})}|_{S_{(\lambda_1,\dots,\lambda_{r-2})}}}{|S_{(\lambda_{r-1},\lambda_r)}|}
    \cdot\N_{S_{(\lambda_1,\dots,\lambda_{r-2})},S_{(\theta_1,\dots,\theta_{r-2})}}(\eta_{(\theta_1,\dots,\theta_{r-2})})\\
&\quad\times \sum_{\substack{(\theta_{r-1},\theta_r)\in\\
\Lambda_{(\lambda_{r-1},\lambda_r)}(\theta'_{r-1})}}
|C_{(\theta_{r-1},\theta_r)}|_{S_{(\lambda_{r-1},\lambda_r)}}\cdot\N_{S_{(\lambda_{r-1},\lambda_r)},S_{(\theta_{r-1},\theta_r)}}(\eta_{(\theta_{r-1},\theta_r)})\\
&=\frac{|S_n|}{|S_{(\lambda_1,\dots,\lambda_{r-2})}|\cdot|S_{(\lambda_{r-1},\lambda_r)}|\cdot|C_\alpha|_{S_n}}\\
&\quad\times\sum_{\theta\in\Lambda_{\lambda}(\alpha)}
 |C_{(\theta_1,\dots,\theta_{r-2})}|_{S_{(\lambda_1,\dots,\lambda_{r-2})}}\cdot
 |C_{(\theta_{r-1},\theta_r)}|_{S_{(\lambda_{r-1},\lambda_r)}}\\
&\quad\times\N_{S_{(\lambda_1,\dots,\lambda_{r-2})},S_{(\theta_1,\dots,\theta_{r-2})}}(\eta_{(\theta_1,\dots,\theta_{r-2})})\cdot
 \N_{S_{(\lambda_{r-1},\lambda_r)},S_{(\theta_{r-1},\theta_r)}}(\eta_{(\theta_{r-1},\theta_r)})\\
&=\quad\frac{|S_n|}{|S_\lambda|\cdot|C_\alpha|_{S_n}}
 \sum_{\theta\in\Lambda_{\lambda}(\alpha)}
 |C_\theta|_{S_\lambda}\cdot\N_{S_\lambda,S_\theta}(\eta_\theta),
\end{align*}
which completes the proof.
\end{proof}

\begin{rem}
  {Note that when $\lambda \ne (n)$, the norms in the summation in the statement of Theorem \ref{thm:proj.balpha.trickle} are only up to a proper subgroup $S_\lambda$ of $S_n$. Therefore, since these multipartitions $\theta\in\Lambda_\lambda(\alpha)$ are not conjugate as $\lambda$--multipartitions, Proposition \ref{lem:norms.from.bips} does not apply and the norms cannot be pulled out of the summation.}
\end{rem}

\begin{thm}\label{thm:coeff.Glambda.in.balpha}
Let $b_\alpha\in B$ and let $\Gamma_\lambda$ be a class element.  Then
\[
b_\alpha=\sum_{\lambda\vdash n}(1_{S_\lambda})^{S_n}(w_\alpha)\xi^{l_\lambda-l_\alpha}\Gamma_\lambda.
\]
\end{thm}

\begin{proof}
Firstly {recall} that the coefficient of $\Gamma_\lambda$ in $b_\alpha$ is exactly
$\l b_\alpha,\T_{w_\lambda}\r$ {by Lemma \ref{lem:class.elt.coeff.in.central.elt}}.  The critical observation is that this value is preserved under the
projection $\pi_\lambda(b_\alpha)$ onto $\H_\lambda${, so that} $\l b_\alpha,\T_{w_\lambda}\r=\l
\pi_\lambda(b_\alpha),\T_{w_\lambda}\r$. By linearity of the inner product we have
\begin{align*}
\l b_\alpha,\T_{w_\lambda}\r
    &=\l \pi_\lambda(b_\alpha),\T_{w_\lambda}\r\\
    &=\frac{|S_n|}{|S_\lambda|\cdot|C_\alpha|_{S_n}}\sum_{\theta\in\Lambda_\lambda(\alpha)}|C_\theta|_{S_\lambda}\l
    \N_{S_\lambda,S_\theta}(\eta_\theta), \T_{w_\lambda}\r.
\end{align*}
Now regardless of which $\theta\in\Lambda_\lambda(\alpha)$ is taken, $l_\theta=l_\alpha$, and so by applying
Theorem \ref{thm:coeff.cox.elt.in.balpha} to each component of $S_\lambda$, the
coefficient of the canonical Coxeter element $\T_{w_\lambda}$ is $\xi^{l_\lambda-l_\alpha}$ in each case.
Thus $\l\N_{S_\lambda,S_\theta}(\eta_\theta), \T_{w_\lambda}\r=\xi^{l_\lambda-l_\alpha}$ is independent of $\theta$, and so may be {pulled} out of the summation.  Then,
\begin{align}
\frac{|S_n|}{|S_\lambda|\cdot|C_\alpha|_{S_n}} \sum_{\theta\in\Lambda_\lambda(\alpha)}|C_\theta|_{S_\lambda}
&=\frac{|S_n|}{|C_\alpha|_{S_n}} \sum_{\theta\in\Lambda_\lambda(\alpha)}\frac{|C_\theta|_{S_\lambda}}{|S_\lambda|}\label{eq:coeff.main.thm}\\
&=|C_{S_n}(w_\alpha)|\sum_{\theta\in\Lambda_\lambda(\alpha)}\frac{1}{|C_{S_\lambda}(w_\theta)|}\notag\\
&=(1_{S_\lambda})^{S_n}(w_\alpha),\notag
\end{align}
which completes the proof.
\end{proof}

\begin{rem}
Utilizing the map $\H\to\H_q$ as described in the Remark in Section \ref{subsec:Hecke.algebras}, this Theorem can ultimately be translated into a corresponding result over $\Z[q,q^{-1}]$ by considering the analogous norms $b_\alpha$ and class elements $\Gamma_\lambda$ over that ring, and replacing the $\xi$ in the statement of the Theorem with $q^{-1}(q-1)$.  See the Remark after Theorem \ref{thm:coeff.cox.elt.in.balpha}.
\end{rem}

\begin{rem}
  Theorem \ref{thm:coeff.Glambda.in.balpha} generalizes the result for the norm of the identity $b_{(1^n)}$ in Theorem \ref{thm:norm.of.1}, effectively making \ref{thm:norm.of.1} redundant.   It also includes as a special case Theorem \ref{thm:coeff.cox.elt.in.balpha}, which gives the coefficient of the Coxeter class element in $b_\alpha$.
\end{rem}

An immediate generalization is the following:
\begin{cor}
  Let $\lambda\vDash n$ and $\mu\le \lambda$.  Then
\[\N_{S_\lambda,S_\mu}(\eta_\mu)=\sum_{\theta\in\Lambda_\lambda}(1_{S_\theta})^{S_\lambda}(w_\mu)\xi^{l_\theta-l_\mu}\Gamma_\theta.\]
\end{cor}

\section{Examples}\label{sec:example}

In this section, we illustrate Theorem \ref{thm:coeff.Glambda.in.balpha} by explicitly giving norms and class elements in terms of the standard basis for $\H$ when $n=3$; we give tables of coefficients for $3\le n\le 5$ which were obtained using Theorem \ref{thm:coeff.Glambda.in.balpha}; and we demonstrate use of the formula to obtain a coefficient in $S_{10}$.

\subsection{Expansions for norms and class elements when $n=3$}
Explicit expressions for norms and class elements for Hecke algebras of types $S_3$ and $S_4$ are listed in \cite{Jon90} (over $\Z[q,q^{-1}]$) and \cite{Fmb} respectively.  We reproduce them for $n=3$ here:

The norms:
\begin{align*}
  b_{(1,1,1)}&=6\T_1+3\xi\T_{s_1}+3\xi\T_{s_2}+\xi^2\T_{s_1s_2}+\xi^2\T_{s_2s_1}+\xi(3+\xi^2)\T_{s_1s_2s_1}\\
            &=6\Gamma_{(1,1,1)}+3\xi\Gamma_{(2,1)}+\xi^2\Gamma_{(3)}\\
  b_{(2,1)}  &=\T_{s_1}+\T_{s_2}+\xi\T_{s_1s_2}+\xi\T_{s_2s_1}+(1+\xi^2)\T_{s_1s_2s_1}\\
            &=\Gamma_{(2,1)}+\xi\Gamma_{(3)}\\
  b_{(3)}    &=\T_{s_1s_2}+\T_{s_2s_1}+\xi\T_{s_1s_2s_1}\\
            &=\Gamma_{(3)},
\intertext{where the $\Gamma_\lambda$ are the following class elements:}
  \Gamma_{(1,1,1)}&=\T_1\\
  \Gamma_{(2,1)}  &=\T_{s_1}+\T_{s_2}+\T_{s_1s_2s_1}\\
  \Gamma_{(3)}    &=\T_{s_1s_2}+\T_{s_2s_1}+\xi\T_{s_1s_2s_1}.
\end{align*}

\subsection{Coefficient tables when $3\le n\le 5$}
\text{}

In the tables below, the row labelled by $\alpha$ gives the coefficients of the class elements in $b_\alpha$.  That is, the entry in position $(\alpha,\lambda)$ gives the coefficient of $\Gamma_\lambda$ in $b_\alpha$.
\vspace{2mm}

\begin{center}
\begin{tabular}{|c||ccc|}
\hline
    &&&\\[-3mm]
$n=3$ &$(1^3)$  &$(2,1)$&$(3)$\\[1mm]
\hline\hline
    &&&\\[-3mm]
$(1^3)$& 6    &$3\xi$ &$\xi^2$ \\[1mm]
$(2,1)$ &       &$1$    &$\xi$   \\[1mm]
$(3)$   &       &       &1\\[1mm]
\hline
\end{tabular}
\quad
\begin{tabular}{|c||ccccc|}
\hline
    &&&&&\\[-3mm]
$n=4$ &$(1^4)$  &$(2,1^2)$&$(3,1)$&$(2^2)$&$(4)$\\[1mm]
\hline\hline
    &&&&&\\[-3mm]
$(1^4)$& 24 &$12\xi$&$4\xi^2$ &$6\xi^2$&$\xi^3$\\[1mm]
$(2,1^2)$ &     &$2$    &$2\xi$   &$2\xi$&$\xi^2$\\[1mm]
$(3,1)$   &     &       &1      &0&$\xi$\\[1mm]
$(2^2)$   &     &       &       &2&$\xi$\\[1mm]
$(4)$&&&&&1\\[1mm]
\hline
\end{tabular}
\end{center}
\vspace{2mm}

\begin{tabular}{|c||ccccccc|}
\hline
    &&&&&&&\\[-3mm]
$n=5$ &$(1^5)$  &$(2,1^3)$&$(3,1^2)$&$(2^2,1)$&$(4,1)$&$(3,2)$&$(5)$\\[1mm]
\hline\hline
    &&&&&&&\\[-3mm]
$(1^5)$& 120    &$60\xi$&$20\xi^2$&$30\xi^2$&$5\xi^3$&$10\xi^3$&$\xi^4$\\[1mm]
$(2,1^3)$ &     &$6$    &$6\xi$   &$6\xi$   &$3\xi^2$&$4\xi^2$&$\xi^3$\\[1mm]
$(3,1^2)$   &     &       &2      &0        &$2\xi$  &$\xi$&$\xi^2$\\[1mm]
$(2^2,1)$   &     &       &       &2        &$\xi$  &$2\xi$&$\xi^2$\\[1mm]
$(4,1)$     &       &   &           &       &1      &0      &$\xi$\\[1mm]
$(3,2)$     &       &       &       &       &       &2      &$\xi$\\[1mm]
$(5)$&&&&&&&1\\[1mm]
\hline
\end{tabular}

\subsection{A coefficient in $S_{10}$}

Suppose we want the coefficient of $\Gamma_{(5,3,2)}$ in $b_{(3,2,2,1,1,1)}$.
Using the form of the coefficients in  Theorem \ref{thm:coeff.Glambda.in.balpha} from the left hand side of equation \eqref{eq:coeff.main.thm}, this is
\[\frac{|S_{10}|\xi^{l_{(5,3,2)}-l_{(3,2,2,1,1,1)}}}{|S_{(5,3,2)}|\cdot|C_{(3,2,2,1,1,1)}|_{S_{10}}}
\sum_{\theta\in\Lambda_{(5,3,2)}({(3,2,2,1,1,1)})}|C_\theta|_{S_{(5,3,2)}}.
\]
We have
\begin{multline*}
\Lambda_{(5,3,2)}\left({(3,2,2,1,1,1)}\right)=\left\{\left((3,2),(2,1),(1,1)\right),\left((3,2),(1,1,1),(2)\right),\right.\\
\left.\left((3,1,1),(2,1),(2)\right),\left((2,2,1),(3),(1,1)\right),\left((2,1,1,1),(3),(2)\right)\right\},
\end{multline*}
and
\begin{align*}
|C_{((3,2),(2,1),(1,1))}|_{S_{(5,3,2)}}&=\frac{5\cdot 4\cdot 3}{3}\cdot\frac{3\cdot 2}{2}&=60\ \\
|C_{((3,2),(1,1,1),(2))}|_{S_{(5,3,2)}}&=\frac{5\cdot 4\cdot 3}{3}&=20\ \\
|C_{((3,1,1),(2,1),(2))}|_{S_{(5,3,2)}}&=\frac{5\cdot 4\cdot 3}{3}\cdot\frac{3\cdot 2}{2}&=60\ \\
|C_{((2,2,1),(3),(1,1))}|_{S_{(5,3,2)}}&=\frac{(5\cdot 4)(3\cdot 2)}{2\cdot 2}\cdot\frac{3\cdot2\cdot1}{3}&=60\ \\
|C_{((2,1,1,1),(3),(2))}|_{S_{(5,3,2)}}&=\frac{5\cdot 4}{2}\cdot\frac{3\cdot2\cdot1}{3}&=20.
\end{align*}
So our coefficient is
\[\frac{10! \xi^{(4+2+1)-(2+1+1)}}{5!3!2!\frac{10\cdot9\cdot8\cdot7\cdot6\cdot5\cdot4}{3\cdot2\cdot2}}
(60+20+60+60+20)=11\xi^3.\]

\end{document}